\documentclass[11pt]{article}
\usepackage[dvips]{graphicx}
\usepackage{textcomp}
\usepackage{amsbsy}
\usepackage{latexsym}
\usepackage[mathscr]{eucal}
\usepackage{amsfonts,amsmath,amsthm}
\usepackage{amssymb}
\usepackage[usenames]{xcolor}

\usepackage{fullpage}

\begin{document}
\bibliographystyle{plain}
\title
{
Data Assimilation and Sampling in Banach spaces
}
\author{ 
Ronald DeVore, Guergana Petrova, and Przemyslaw Wojtaszczyk
\thanks{%
   This research was supported by the ONR Contracts
   N00014-15-1-2181 and N00014-16-1-2706; 
  the  NSF Grant
       DMS 15-21067;  DARPA Grant  HR0011619523 through Oak Ridge National Laboratory,
    and the  Polish NCN grant DEC2011/03/B/ST1/04902. }  }
\hbadness=10000
\vbadness=10000
\newtheorem{lemma}{Lemma}[section]
\newtheorem{prop}[lemma]{Proposition}
\newtheorem{cor}[lemma]{Corollary}
\newtheorem{theorem}[lemma]{Theorem}
\newtheorem{remark}[lemma]{Remark}
\newtheorem{example}[lemma]{Example}
\newtheorem{definition}[lemma]{Definition}
\newtheorem{proper}[lemma]{Properties}
\newtheorem{assumption}[lemma]{Assumption}
%
\newenvironment{disarray}{\everymath{\displaystyle\everymath{}}\array}{\endarray}

\def\RR{\rm \hbox{I\kern-.2em\hbox{R}}}
\def\NN{\rm \hbox{I\kern-.2em\hbox{N}}}
\def\ZZ{\rm {{\rm Z}\kern-.28em{\rm Z}}}
\def\CC{\rm \hbox{C\kern -.5em {\raise .32ex \hbox{$\scriptscriptstyle
|$}}\kern
-.22em{\raise .6ex \hbox{$\scriptscriptstyle |$}}\kern .4em}}
\def\vp{\varphi}
\def\<{\langle}
\def\>{\rangle}
\def\t{\tilde}
\def\i{\infty}
\def\e{\varepsilon}
\def\sm{\setminus}
\def\nl{\newline}
\def\o{\overline}
\def\wt{\widetilde}
\def\wh{\widehat}
\def\cT{{\cal T}}
\def\cA{{\cal A}}
\def\cI{{\cal I}}
\def\cV{{\cal V}}
\def\cB{{\cal B}}
\def\cF{{\cal F}}
\def\cY{{\cal Y}}

\def\cD{{\cal D}}
\def\cP{{\cal P}}
\def\cJ{{\cal J}}
\def\cM{{\cal M}}
\def\cO{{\cal O}}
\def\Chi{\raise .3ex
\hbox{\large $\chi$}} \def\vp{\varphi}
\def\lsima{\hbox{\kern -.6em\raisebox{-1ex}{$~\stackrel{\textstyle<}{\sim}~$}}\kern -.4em}
\def\lsim{\hbox{\kern -.2em\raisebox{-1ex}{$~\stackrel{\textstyle<}{\sim}~$}}\kern -.2em}
\def\[{\Bigl [}
\def\]{\Bigr ]}
\def\({\Bigl (}
\def\){\Bigr )}
\def\[{\Bigl [}
\def\]{\Bigr ]}
\def\({\Bigl (}
\def\){\Bigr )}
\def\L{\pounds}
\def\pr{{\rm Prob}}
\newcommand{\cs}[1]{{\color{magenta}{#1}}}
\def\ds{\displaystyle}
\def\ev#1{\vec{#1}}     
\newcommand{\lt}{\ell^{2}(\nabla)}
\def\Supp#1{{\rm supp\,}{#1}}
\def\R{\mathbb{R}}
\def\E{\mathbb{E}}
\def\nl{\newline}
\def\T{{\relax\ifmmode I\!\!\hspace{-1pt}T\else$I\!\!\hspace{-1pt}T$\fi}}
\def\N{\mathbb{N}}
\def\Z{\mathbb{Z}}
\def\N{\mathbb{N}}
\def\Zd{\Z^d}
\def\Q{\mathbb{Q}}
\def\C{\mathbb{C}}
\def\Rd{\R^d}
\def\gsim{\mathrel{\raisebox{-4pt}{$\stackrel{\textstyle>}{\sim}$}}}
\def\sime{\raisebox{0ex}{$~\stackrel{\textstyle\sim}{=}~$}}
\def\lsim{\raisebox{-1ex}{$~\stackrel{\textstyle<}{\sim}~$}}
\def\div{\mbox{ div }}
\def\M{M}  \def\NN{N}                  
\def\L{{\ell}}               
\def\Le{{\ell^1}}            
\def\Lz{{\ell^2}}
\def\Let{{\tilde\ell^1}}     
\def\Lzt{{\tilde\ell^2}}
\def\Ltw{\ell^\tau^w(\nabla)}
\def\t#1{\tilde{#1}}
\def\la{\lambda}
\def\La{\Lambda}
\def\ga{\gamma}
\def\BV{{\rm BV}}
\def\Ga{\eta}
\def\al{\alpha}
\def\cZ{{\cal Z}}
\def\cA{{\cal A}}
\def\cU{{\cal U}}
\def\argmin{\mathop{\rm argmin}}
\def\argmax{\mathop{\rm argmax}}
\def\prob{\mathop{\rm prob}}

\def\cO{{\cal O}}
\def\cA{{\cal A}}
\def\cC{{\cal C}}
\def\cS{{\cal F}}
\def\bu{{\bf u}}
\def\bz{{\bf z}}
\def\bZ{{\bf Z}}
\def\bI{{\bf I}}
\def\cE{{\cal E}}
\def\cD{{\cal D}}
\def\cG{{\cal G}}
\def\cI{{\cal I}}
\def\cJ{{\cal J}}
\def\cM{{\cal M}}
\def\cN{{\cal N}}
\def\cT{{\cal T}}
\def\cU{{\cal U}}
\def\cV{{\cal V}}
\def\cW{{\cal W}}
\def\cL{{\cal L}}
\def\cB{{\cal B}}
\def\cG{{\cal G}}
\def\cK{{\cal K}}
\def\cX{{\cal X}}
\def\cS{{\cal S}}
\def\cP{{\cal P}}
\def\cQ{{\cal Q}}
\def\cR{{\cal R}}
\def\cU{{\cal U}}
\def\bL{{\bf L}}
\def\bl{{\bf l}}
\def\bK{{\bf K}}
\def\bC{{\bf C}}
\def\X{X\in\{L,R\}}
\def\ph{{\varphi}}
\def\D{{\Delta}}
\def\H{{\cal H}}
\def\bM{{\bf M}}
\def\bx{{\bf x}}
\def\bj{{\bf j}}
\def\bG{{\bf G}}
\def\bP{{\bf P}}
\def\bW{{\bf W}}
\def\bT{{\bf T}}
\def\bV{{\bf V}}
\def\bv{{\bf v}}
\def\bt{{\bf t}}
\def\bz{{\bf z}}
\def\bw{{\bf w}}
\def \span{{\rm span}}
\def \meas {{\rm meas}}
\def\rhom{{\rho^m}}
\def\diff{\hbox{\tiny $\Delta$}}
\def\EE{{\rm Exp}}
\def\lll{\langle}
\def\argmin{\mathop{\rm argmin}}
\def\codim{\mathop{\rm codim}}
\def\rank{\mathop{\rm rank}}

\def\argmax{\mathop{\rm argmax}}
\def\dJ{\nabla}
\newcommand{\ba}{{\bf a}}
\newcommand{\bb}{{\bf b}}
\newcommand{\bc}{{\bf c}}
\newcommand{\bd}{{\bf d}}
\newcommand{\bs}{{\bf s}}
\newcommand{\bff}{{\bf f}}
\newcommand{\bp}{{\bf p}}
\newcommand{\bg}{{\bf g}}
\newcommand{\by}{{\bf y}}
\newcommand{\br}{{\bf r}}
\newcommand{\be}{\begin{equation}}
\newcommand{\ee}{\end{equation}}
\newcommand{\bea}{$$ \begin{array}{lll}}
\newcommand{\eea}{\end{array} $$}
\def \Vol{\mathop{\rm  Vol}}
\def \mes{\mathop{\rm mes}}
\def \Prob{\mathop{\rm  Prob}}
\def \exp{\mathop{\rm    exp}}
\def \sign{\mathop{\rm   sign}}
\def \sp{\mathop{\rm   span}}
\def \rad{\mathop{\rm   rad}}
\def \vphi{{\varphi}}
\def \csp{\overline \mathop{\rm   span}}
%
%
\newcommand{\beqn}{\begin{equation}}
\newcommand{\eeqn}{\end{equation}}
\def\beginproof{\noindent{\bf Proof:}~ }
\def\endproof{\hfill\rule{1.5mm}{1.5mm}\\[2mm]}

\newenvironment{Proof}{\noindent{\bf Proof:}\quad}{\endproof}

\renewcommand{\theequation}{\thesection.\arabic{equation}}
\renewcommand{\thefigure}{\thesection.\arabic{figure}}

\makeatletter
\@addtoreset{equation}{section}
\makeatother

\newcommand\abs[1]{\left|#1\right|}
\newcommand\clos{\mathop{\rm clos}\nolimits}
\newcommand\trunc{\mathop{\rm trunc}\nolimits}
\renewcommand\d{d}
\newcommand\dd{d}
\newcommand\diag{\mathop{\rm diag}}
\newcommand\dist{\mathop{\rm dist}}
\newcommand\diam{\mathop{\rm diam}}
\newcommand\cond{\mathop{\rm cond}\nolimits}
\newcommand\eref[1]{{\rm (\ref{#1})}}
\newcommand{\iref}[1]{{\rm (\ref{#1})}}
\newcommand\Hnorm[1]{\norm{#1}_{H^s([0,1])}}
\def\int{\intop\limits}
\renewcommand\labelenumi{(\roman{enumi})}
\newcommand\lnorm[1]{\norm{#1}_{\ell^2(\Z)}}
\newcommand\Lnorm[1]{\norm{#1}_{L_2([0,1])}}
\newcommand\LR{{L_2(\R)}}
\newcommand\LRnorm[1]{\norm{#1}_\LR}
\newcommand\Matrix[2]{\hphantom{#1}_#2#1}
\newcommand\norm[1]{\left\|#1\right\|}
\newcommand\ogauss[1]{\left\lceil#1\right\rceil}
\newcommand{\QED}{\hfill
\raisebox{-2pt}{\rule{5.6pt}{8pt}\rule{4pt}{0pt}}%
  \smallskip\par}
\newcommand\Rscalar[1]{\scalar{#1}_\R}
\newcommand\scalar[1]{\left(#1\right)}
\newcommand\Scalar[1]{\scalar{#1}_{[0,1]}}
\newcommand\Span{\mathop{\rm span}}
\newcommand\supp{\mathop{\rm supp}}
\newcommand\ugauss[1]{\left\lfloor#1\right\rfloor}
\newcommand\with{\, : \,}
\newcommand\Null{{\bf 0}}
\newcommand\bA{{\bf A}}
\newcommand\bB{{\bf B}}
\newcommand\bR{{\bf R}}
\newcommand\bD{{\bf D}}
\newcommand\bE{{\bf E}}
\newcommand\bF{{\bf F}}
\newcommand\bH{{\bf H}}
\newcommand\bU{{\bf U}}
\newcommand\cH{{\cal H}}
\newcommand\sinc{{\rm sinc}}
\def\enorm#1{| \! | \! | #1 | \! | \! |}

\newcommand{\dm}{\frac{d-1}{d}}

\let\bm\bf
\newcommand{\bbeta}{{\mbox{\boldmath$\beta$}}}
\newcommand{\bal}{{\mbox{\boldmath$\alpha$}}}
\newcommand{\bbi}{{\bm i}}

\def\nnew{\color{Red}}
\def\mnew{\color{Blue}}

\newcommand{\dI}{\Delta}
\newcommand\aconv{\mathop{\rm absconv}}

\maketitle
\date{}
\begin{abstract}

This paper studies the problem of approximating a function $f$ in a Banach space $\cX$ from measurements $l_j(f)$, $j=1,\dots,m$, where the $l_j$
are linear functionals from $\cX^*$.
  Quantitative results  for such recovery problems require additional information about the sought after function $f$.  These additional assumptions take the form of
assuming that $f$ is in a certain model class $K\subset \cX$.
      Since there are generally infinitely many functions in $K$ which share 
these same measurements, the best approximation is the center of the smallest ball $B$, called the {\it Chebyshev ball}, which contains the set 
$\bar K$ of all  $f$ in $K$ with these measurements.  Therefore, the problem is reduced to analytically or numerically approximating this Chebyshev ball.

Most  results    study this problem for classical Banach spaces $\cX$ such as the $L_p$ spaces, 
$1\le p\le \infty$,  and  for $K$   the unit ball
of a smoothness space in $\cX$.    Our interest in this paper is in the   model classes $K=\cK(\e,V)$, with  $\e>0$ and $V$   a finite dimensional subspace of $\cX$,  which consists of all $f\in \cX$ such that
$\dist(f,V)_\cX\le \e$.   These model classes, called {\it approximation sets},  arise naturally in application domains such as parametric partial differential equations, uncertainty quantification, and signal processing.

A  general theory for the   recovery of approximation sets  in a Banach space  is given.  This theory includes tight a priori bounds on optimal performance, and
algorithms for finding near optimal approximations.      It  builds on the initial 
 analysis given in  \cite{MPPY}
for the case when  $\cX$ is  a Hilbert space,  and further studied in \cite{BCDDPW1}.   
 It is shown how the recovery problem for approximation sets is  connected with  well-studied concepts in Banach space theory such as liftings and the angle between spaces. 
Examples are given that show how this theory can be used to recover several recent results on sampling and data assimilation.
  \end{abstract}

\section{Introduction}
\label{Introduction}

 One of the most ubiquitous problems in science is to approximate an unknown function $f$ from given data observations of  $f$.
 Problems of this type go under the  terminology of {\it optimal recovery},  {\it data assimilation} or the more colloquial terminology of {\it data fitting}.   
 To prove quantitative results  about the accuracy of such recovery requires, not only the data observations, but also additional information about the sought after function $f$.  Such additional information takes the form of
assuming that $f$ is in a prescribed model class $K\subset \cX$.

     The classical setting for such problems    (see, for example, \cite{BB,MR,MRW, TW1})
is that one has a  bounded set $K$ in a Banach space $\cX$ and
a finite collection of linear functionals $l_j$, $j=1,\dots,m$, from $\cX^*$.   Given a function which is known to be in $K$
and to have known measurements $l_j(f)=w_j$, $j=1,\dots,m$, the {\it optimal recovery } problem is to construct the best
approximation to $f$ from this information.   Since there are generally infinitely many functions in $K$ which share 
these same measurements, the best approximation is the center of the smallest ball $B$, called the {\it Chebyshev ball}, which contains the set 
$\bar K$ of all  $f$ in $K$ with these measurements.  The best error of approximation is then the radius of this Chebyshev ball.

Most  results in optimal recovery  study this problem for classical Banach spaces $\cX$ such as the $L_p$ spaces, 
$1\le p\le \infty$,  and  for $K$   the unit ball
of a smoothness space in $\cX$.    Our interest in this paper is in certain other model classes $K$, called {\it approximation sets}, that arise in various applications.
As a motivating example, consider the analysis of  complex physical systems from data observations.  In such settings the sought after functions    satisfy a (system of) parametric partial differential
equation(s) with unknown parameters and hence lie in the solution manifold $\cM$ of the parametric model.  There may also be uncertainty in the parametric model.   Problems of this type fall into the general paradigm of uncertainty quantification.
The solution manifold of a parametric partial differential equation (pde) is a complex object and information about the manifold is usually only known through approximation results on how well  the  elements in the manifold can be approximated by certain low dimensional linear spaces or by  sparse Taylor (or other)  polynomial expansions (see \cite{CD}).   For this reason, the manifold $\cM$ is often  replaced, as a model class, by 
the set
$\cK(\e,V)$ consisting of  all elements in $\cX$  that can be   approximated by the  linear space $V=V_n$ of   dimension $n$ to  
accuracy $\e=\e_n$.    We call these model classes $\cK(\e,V)$  {\it approximation sets} and they are the main focus of this paper.  Approximation sets also arise naturally as model classes in other settings such as signal processing where the problem is to construct an approximation to a signal from samples (see e.g. \cite{AHP, AH,AHS, DT, TW} and the papers cited therein
as representative), although this terminology is not in common use in that setting.

  Optimal  recovery in this new setting of approximation sets as the model class was formulated 
and analyzed in \cite{MPPY}
when $\cX$ is  a Hilbert space,  and further studied in \cite{BCDDPW1}.  In particular,  it was shown in the latter paper 
that a certain numerical algorithm proposed
in \cite{MPPY}, based on least squares approximation,   is optimal.

The purpose of the present paper is to  provide a general theory for the optimal  or near optimal recovery of    approximation sets in  a general  Banach space $\cX$. 
   While, as noted in the abstract,  the optimal recovery has a simple theoretical description as the center of
the Chebyshev ball
and the optimal performance, i.e., the best error,  is given by the radius of the Chebyshev ball, this is far from a 
satisfactory solution to the problem since it is not clear how to find the  center and the radius of  the Chebyshev ball.  
This leads to the two fundamental problems studied in the  paper.
The first centers on building numerically executable   algorithms which are optimal or perhaps  only near optimal, i.e., they either determine the Chebyshev center or approximate it
sufficiently well. 
The second  
problem is to give sharp a priori bounds for the best error, 
 i.e the Chebyshev radius, 
in terms of easily computed 
quantities.  We show how these two problems are connected with well-studied concepts in Banach space theory such as liftings and the angle between spaces. 
Our main results determine a priori bounds for optimal algorithms and give numerical recipes for obtaining optimal or near optimal algorithms.   

\subsection{Problem formulation and summary of results}

Let $\cX$  be a Banach space with norm $\|\cdot\|=\|\cdot\|_\cX$ and let  $S\subset \cX$ be  any  subset of $\cX$. 
We assume we are   given  measurement functionals $l_1,\dots,l_m\in \cX^*$ that   are linearly independent.   
We study the general question of how best to recover a function $f$ from the information that $f\in S$ and $f$ 
has the known measurements 
$M(f):= M_m(f):=(l_1(f), \ldots,l_m(f))=(w_1,\ldots,w_m)\in \R^m$.

An algorithm $A$
for this  recovery problem is a mapping which when given the measurement data  $w=(w_1,\ldots,w_m)$
assigns an element $ A(w)\in \cX$ as the approximation to $f$.  Thus, an algorithm is a possibly nonlinear mapping
$$
A: \  \R^m \mapsto \cX.
$$
Note that there are generally  many functions $f\in S$ which share the same data.  We denote this collection by
$$
S_w:= \{f\in S:\ l_j(f)=w_j,\ j=1,\dots,m\}.
$$
A {\it pointwise optimal algorithm }$  A^*$  (if it exists) is one which minimizes the worst error for each $w$:
$$
A^*(w):=\argmin_{g\in\cX}\sup_{f\in S_w}\|f-g\|.
$$
This optimal algorithm has a simple geometrical description that is well known (see e.g. \cite{MR}).  For a given $w$, we consider all balls   $B(a,r)\subset \cX$  
which contain $S_w$.  The smallest     ball  $B(a(w),r(w))$, if it exists, is called the  {\it Chebyshev ball} and its radius is called the {\it Chebyshev radius} of $S_w$.
We postpone the discussion of existence, uniqueness, and properties of the smallest ball to the next section.  
For now, we remark
that when the Chebyshev   ball  $B(a(w),r(w))$ exists for each measurement vector $w\in\R^m$, then   the
 pointwise optimal 
algorithm is 
the mapping $A^*:\ w\to a(w)$ and the      pointwise optimal error   for this recovery problem is given by 
\be
\label{besterror}
 \rad(S_w):=\sup_{f\in S_w}\|f-a(w)\|=  \inf_{a\in\cX} \sup_{f\in S_w}\|f-a\|= \inf_{a\in\cX} \inf_r\{S_w\subset B(a,r)\}.
\ee
 In summary, the smallest error that any algorithm for the recovery of $S_w$
can attain is $\rad(S_w)$, and it is attained by taking the center of the Chebyshev ball of $S_w$. 
\vskip .1in 
\noindent
{\bf Pointwise Near  Optimal Algorithm:}  {\it We say that an algorithm $A$ is pointwise near  optimal with constant $C$ for the  set  $S$  if 
$$
\sup_{f\in S_w}\|f-A(w)\|\le C\rad(S_w),  \quad \forall w\in\R^m.
$$
}
\vskip .1in
 
In applications, one typically knows the linear functionals $l_j$, $j=1,\dots,m$, (the measurement devices) but has no a priori 
knowledge of  the measurement values $w_j$, $j=1,\dots,m$, that will arise.   Therefore, a second meaningful measure of performance is
\be
\label{defR}
R(S):=\sup_{w\in\R^m}\rad(S_w).
\ee
Note that an algorithm which achieves the bound $R(S)$  will generally not be optimal for each $w\in \R^m$.
  We refer to the second type of estimates as global and  a global  optimal    algorithm would be 
one that achieved the bound $R(S)$. 
\vskip .1in
\noindent
{\bf Global Near  Optimal   Algorithm:}  {\it We say that an algorithm $A$ is  a global  near  optimal   algorithm with constant $C$ for the set $S$,  if 
$$
\sup_{w\in \R^m} \sup_{f\in S_w}\|f-A(w)\|\le C R(S).
$$
}
\vskip .1in
 \noindent
  Note that if an algorithm is near  optimal  for each of the sets  $S_w$ with a constant $C$,  independent of $w$, then it is a  global near  optimal    algorithm with the same constant $C$.

The above description in terms of $\rad(S_w)$ provides a nice simple geometrical description of the optimal recovery problem. 
However, it is not a practical solution for a given set $S$, since the problem of finding the Chebyshev center  and radius of $S_w$ is essentially the same as the original optimal recovery problem, and moreover, is known, in general,  to  be 
NP hard (see \cite{NP}).  Nevertheless, it provides some guide to the construction of optimal or near optimal algorithms.

 In the first part of this paper, namely \S\ref{sec:prelim} and \S\ref{sec:liftings},  we use classical ideas and techniques of Banach space theory 
 (see, for example,  \cite{LT,BL,Singer, Wbook}), to provide results on the  optimal recovery of the sets $S_w$ for any set $S\subset \cX$,  $w\in\R^m$.   Much of the material
 in these two sections is known or easily derived from known results, but we recount this for the benefit of the reader and to ease the exposition of this paper.
 
   Not surprisingly,
the form of these results depends very much on the structure of the Banach space.   Let us recall that the unit ball $U:=U(\cX)$  of the 
Banach space $\cX$ is always convex.   The Banach space $\cX$ is said to be {\it strictly  convex} if  
$$
\|f+g\|< 2, \quad  \forall f,g\in U, \quad  f\neq g.
$$
A stronger property of $\cX$  is the {\it uniform convexity}. To describe this property, we introduce the modulus of convexity of $\cX$  defined  by
\be
\label{uc}
\delta(\e):=\delta_\cX(\e):=\inf  
  \left\{ 1 - \left\| \frac{f + g}{2} \right\| \,:\,  f, g \in U, \| f - g \| \geq \e \right\},\quad \e >0.
\ee
The space $\cX$ is  called uniformly convex if $\delta(\e)>0$ for all $\e>0$.  For uniformly convex spaces $\cX$, it  is known that $\delta$ is a strictly  increasing function taking values in $[0,1]$ as $\e$ varies in $[0,2] $  (see, for example,   \cite[Th. 2.3.7]{AOS}). 

Uniform convexity implies strict convexity and it also implies that $\cX$ is reflexive (see [Prop. 1.e.3] in \cite{LT}), i.e. $\cX^{\ast\ast}=\cX$, where $\cX^*$ denotes the dual
space of $\cX$.  If $\cX$ is uniformly convex, then there is quite a similarity between the results we obtain and those in the Hilbert space case.  This covers, for example, the case when 
$\cX$ is an $L_p$ or $\ell_p$ space for $1<p<\infty$, or one of the Sobolev spaces  for these $L_p$.  The cases $p=1,\infty$, as well as the case of a general Banach space $\cX$, add some new wrinkles and the theory is not as complete.

The next  part of this paper turns to  the model classes of main interest to us:
\vskip .1in

\noindent
{\bf Approximation Set:} 
{\it  We call the set $\cK=\cK(\e,V)$ an approximation set  if
$$
\cK=\{f\in \cX\,:\,\,\dist(f, V)_\cX\le \e\},
$$
where $V\subset \cX$ is a known finite dimensional  space.
}
\vskip .1in
\noindent
We denote the dependence of $\cK$ on $\e$ and $V$ only when this is not clear from the context.

The main contribution  of the present paper is to  describe near optimal  algorithms for the recovery of approximation sets in a   general Banach space $\cX$.  The determination of  optimal  or near  optimal  algorithms  and their performance for approximation sets is connected to
  liftings (see \S \ref{sec:liftings}),   and the angle between the space $V$ and the null space $\cN$ of the measurements (see \S \ref{sec:radius}).
These concepts allow  us to describe a general procedure  for constructing   recovery algorithms $A$ which are  pointwise near  optimal  (with constant $2$)   and hence are  also globally near  optimal   (see \S \ref{sec:algorithms}).    The near optimal algorithms  $A: \R^m\mapsto \cX$ that we construct satisfy the performance bound
\be
\label{nopt}
\|f-A(M(f))\|_\cX\le C \mu(\cN,V)\e,\quad f\in\cK(\e,V),
\ee
where $\mu(V,\cN)$ is the reciprocal of the angle $\theta(\cN,V)$ between $V$ and the null space $\cN$ of the measurement map $M$ and $C$ is any constant larger than $4$.   We prove that this estimate is near optimal in the sense that, for any recovery algorithm $A$, we have
\be
\label{nopt1}
\sup_{f\in \cK(\e,V)} \|f-A(M(f))\|_\cX\ge  \mu(cN,V)\e, 
\ee
and so the only possible improvement that can be made in our algorithm is to reduce the size of the constant $C$.  Thus, the constant $\mu(V,\cN)$, respectively the angle $\theta(V,\cN)$
determines how well we can recover approximation sets and quantifies the compatibility between    the measurements and $V$.    As we have already noted, this result is not new for the Hilbert space setting.   We discuss in \S \ref{sec:sampling}, other settings where this angle has been recognized to determine best recovery.
 
 It turns out that the reconstruction algorithm $A$ we provide, does not depend on $\e$ (and, in fact,  does not require the knowledge of $\e$) and therefore gives the bound
$$
\|f-A(M(f))\|_\cX \le C \mu(\cN,V) \dist(f,V)_\cX,\quad f\in \cX.
$$

 In section \S \ref{sec:performance}, we  give examples of how to implement 
 our  near  optimal  algorithm in concrete settings  when $\cX=L_p$, or $\cX$ is the space of continuous functions on a domain $D\subset \R^d$, $\cX=C(D)$.  As a representative example of the results in that section, consider  the case $\cX=C(D)$ with $D$ a domain in $\R^d$, and suppose that the  measurements functionals are 
  $l_j(f)=f(P_j)$, $j=1,\dots,m$, where the $P_j$ are points in $D$.  Then, we prove that
$$
\frac{1}{2} \mu(\cN,V) \le   \sup_{v\in V}\frac{\|v\|_{C(D)}}{\displaystyle{\max_{1\le j\le m}|v(P_j)}|} \le 2 \mu(\cN,V).
$$
Hence, the performance of this data fitting problem is controlled by the ratio of the continuous and discrete norms on $V$.  
 Results of this type are well-known, via Lebesgue constants,
in the case of interpolation (when $m=n$).   

In \S \ref{sec:sampling}, we discuss how our results are related to generalized sampling in a Banach space and discuss how several recent results in sampling can be obtained from our 
approach.   Finally, in \S \ref{sec:choosemeasurements},   we discuss good choices for where to take measurements if this option is available to us.


\section{Preliminary remarks }
\label{sec:prelim}   
 
In this section, we recall some standard concepts in Banach spaces and relate them to the optimal recovery problem of interest to us.  We work in the setting that $S$ is any set (not necessarily an approximation set).   The results of this section are essentially known
and are given only to orient the reader.
 \subsection{The Chebyshev ball} 
 \label{ss:Chebrad}

   Note that, in general,   the center of the Chebyshev ball may not come from $S$.   This can even occur in finite dimensional setting (see  
   Example \ref{Exnonsymmetric} given below).    However, it may be desired, or even required in certain applications,  that the recovery for $S$ 
   be a point from $S$.  The description of optimal  algorithms with this requirement is connected with what we call the  {\it restricted Chebyshev ball of $S$}. 
   To explain   this, we introduce some further
   geometrical concepts.

 For a   subset $S$ in a Banach space $\cX$,  we define the following quantities:
  \begin{itemize}
 \item The {\bf diameter } of $S$ is defined by $\diam(S):= \displaystyle{\sup_{f,g\in S}\|f-g\|}$.
\item The {\bf restricted  Chebyshev radius } of $S$  is defined by 
$$
{\rad}_C(S):=\inf_{a\in S}\inf_{r}\{S \subset B(a,r)\}=\inf_{a\in S}\sup_{f\in S} \|f-a\|.
$$
\item  The  {\bf  Chebyshev radius} of $S$  was already defined as
$$
\rad(S):=\inf_{a\in \cX}\inf_{r}\{S \subset B(a,r)\}=\inf_{a\in \cX}\sup_{f\in S} \|f-a\|.
$$
\end{itemize}

\noindent
\begin{remark}
\label{tri}
It is clear that for every $S\subset \cX$,  we have 
\begin{equation}\label{A1}
\diam (S)\geq {\rad}_C(S)\geq \rad(S)\geq \tfrac12 \diam(S).
\end{equation}
\end{remark}
\noindent
Let us start with the following theorem, that tells us that we can construct near  optimal algorithms for the recovery of the  set $S$ if we can simply find a point $a\in S$.
   \begin{theorem}
   \label{theorem:nearoptimal}
   Let $S$ be any subset of $\cX$.   If $a\in S$, then
   \be
   \label{no}
 {\rad} (S)\le  {\rad}_C(S)\leq \sup_{f\in S} \|f-a\|\le 2 \rad(S), 
      \ee
   and therefore $a$ is, up to the constant $2$, an optimal recovery of $S$.  
  \end{theorem}
  
  \noindent
  {\bf Proof:}   Let $B(a',r)$, $a'\in\cX$, be any ball that contains $S$.   Then, for any $f\in S$,
 $$
  \|f-a\|\le \|f-a'\|+\|a-a'\|\le 2r.
 $$
  Taking an infimum over all such balls we obtain the theorem.\hfill $\blacksquare$

   We say that an $a\in \cX$ which recovers $S$ with the accuracy of \eref{no}  provides  a {\it near optimal} recovery with constant $2$.
   We shall use this theorem in our construction of algorithms for the recovery of the sets $\cK_w $, when $\cK$ is an approximation set.
 The relevance of the above theorem and remark  viz a viz for our recovery problem is that if we determine the diameter or restricted  Chebyshev radius of  $\cK_w$, we will determine the optimal
error $\rad(\cK_w)$ in the recovery problem,  but only up to the factor two.     

     
 \subsection{Is $\rad(S)$ assumed?} 
 \label{ss:assumed}

  In view of  the discussion preceding \eref{besterror}, the best  pointwise error we can achieve by any recovery algorithm for $S$ is given by  $\rad(S)$.   Unfortunately,  in general, for arbitrary bounded sets $S$ in a general infinite dimensional Banach space $\cX$, the radius $\rad(S)$ may not be assumed.
  The first such example was given  in \cite{G}, where    $\cX=\{f\in C[-1,1] :\int_{-1}^1 f(t)\, dt =0\}$ with the uniform norm and    $S\subset \cX$ is a set consisting of three functions.  Another, 
  simpler example of a set $S$ in this same space was given in \cite{SWard}.  In \cite{K},  it is shown that each nonreflexive space admits an equivalent norm for which such examples also exist.
  If we place more structure on the Banach space $\cX$, then we can show that the radius of any bounded subset $S\subset \cX$ is assumed.    We present the following special case of an old result of Garkavi (see  \cite[Th. II]{G}).
  
\begin{lemma}
\label{lemma:assumed}   If the Banach space $\cX$ is reflexive  (in particular,  if it is finite dimensional),  then for any bounded set $S\subset \cX$, $\rad(S)$ is assumed in the sense that there is a   ball $B(a,r)$ with $r=\rad(S)$ which contains $S$.  If, in addition,  $\cX$ is uniformly convex,  then this ball is unique.
 \end{lemma} 
 \noindent
{\bf Proof:}    Let $B(a_n,r_n)$ be balls which contain $S$ and for which $r_n\to \rad(S)=:r$. 
Since $S$ is bounded,  the $a_n$ are bounded and hence, without loss of generality, 
 we can assume that $a_n$ converges weakly to $a\in \cX$
 (since every bounded sequence in a reflexive Banach space has a weakly converging subsequence).  Now let $f$ be any element in $S$.
Then, there is a norming functional $l\in \cX^*$ of norm one for which $l(f-a)=\|f-a\|$.   Therefore
$$
\|f-a\|=l(f-a)=\lim_{n\to\infty} l(f-a_n)\le \lim_{n\to\infty}\|f-a_n\|\le \lim_{n\to\infty } r_n=r.
$$
This shows that $B(a,r)$ contains $S$ and so the radius is attained.  If  $\cX$ is uniformly convex  and we assume that there are two balls, centered at $a$ and $a'$,
 respectively, $a\neq a'$, each of radius $r$ which contain $S$.  
If   $\e:=\|a-a'\|>0$, since $\cX$ is  uniformly convex, for every $f\in S$ and   for $\bar a:=\frac{1}{2}(a+a')$, we have
$$\|f-\bar a\|=\left\|\frac{f-a}{2}+\frac{f-a'}{2}\right\|\le r- r\delta(\e/r)<r,
$$
which contradicts the fact that  $\rad(S)=r$. \hfill $\blacksquare$

\subsection{Some examples}
\label{exmp}  
 In this section, we will show that   for centrally symmetric, convex sets $S$, we have a very explicit relationship between the
 $\diam(S)$, $\rad(S)$, and ${\rad}_C(S)$. We also  give some examples showing  that for general sets  $S$ the situation is more involved and the only relationship between the latter quantities is the one given by Remark \ref{tri}.

\begin{prop}
\label{prop:centrallysymmetric}
Let $S\subset \cX$ be a centrally symmetric, convex set in a Banach space $\cX$.  Then, we have
 
\noindent
{\rm (i)}  the smallest ball containing $S$ is centered at $0$ and 
$$
\diam(S)=2\sup_{f\in S} \|f\| =2 \rad(S)=2 {\rad}_C(S).
$$

\noindent
{\rm (ii)}  for any $w\in \R^m$,  we have
$\diam (S_w)\leq \diam (S_0)$.
 \end{prop}
 \noindent
{\bf Proof:}  (i) We need only consider the case  when $r:=\sup_{f\in S} \|f\|<\infty$. 
 Clearly $0\in S$, $S\subset B(0,r)$, and thus ${\rad}_C(S)\le  r$.
In  view of \eref{A1},
$\diam(S)\le 2\rad(S)\le 2{\rad}_C(S)\le  2r$.
Therefore, we need only show that $\diam(S)\ge 2r$.   For  any  $\e>0$,  let $f_\e\in \cS$ be such that $\|f_\e\|\geq r-\e$.  
Since $S$ is centrally symmetric $-f_\e\in S$ and 
$\diam(S)\geq \|f_\e-(-f_\e)\|\geq 2r-2\e.$ Since $\e>0$ is arbitrary, $\diam(S)\ge 2r$, as desired. 
  \newline (ii) We need only consider the case $\diam(S_0)<\infty$.   Let $a,b\in S_w$.    
  From the convexity and central symmetry of $S$,  we know that $\frac{1}{2}(a-b)$ and  
  $\frac{1}{2}(b-a)$ are both  in $S_0$.   Therefore
 $$
 \diam S_0\geq  \|\frac12(a-b)-\frac12(b-a)\|=\|a-b\|.
$$
  Since $a,b$ were arbitrary,  we get $\diam(S_0)\geq\diam (S_w)$.
   \hfill $\blacksquare$

   In what follows, we denote by $\ell_p(\N)$ the set of all real sequences $x$, such that 
$$
\ell_p(\N):=\{x=(x_1,x_2,\ldots):\,\|x\|_{\ell_p(\N)}:=(\sum_{j=1}^\infty|x_j|^p)^{1/p}<\infty\}, \quad 1\leq p<\infty,
$$
$$
\ell_\infty(\N):=\{x=(x_1,x_2,\ldots):\,\|x\|_{\ell_\infty(\N)}:=\sup_j|x_j|<\infty\},
$$
and 
$$
c_0:=\{x=(x_1,x_2,\ldots):\lim_{j\rightarrow 0} x_j=0\}, \quad \|x\|_{c_0}=\|x\|_{\ell_\infty(\N)}.
$$

\noindent
We start with the following example which can be found  in \cite{BKS}.
\begin{example}
Let  $\cX =\ell_2(\N)$ with a  new  norm defined as
$$
\|x\|=\max\{\tfrac 12 \|x\|_{\ell_2(\N)},\|x\|_{\ell_\infty(\N)}\},\quad x\in \ell_2(\N),
$$
and consider the set $S:=\{x \ :\ \|x\|_{\ell_2(\N)}\leq 1 \ \mbox{ and }\  x_j\geq 0 \ \mbox{ for all } \ j\}$.  Then, for this set $S$, 
 $$\diam(S)=\rad(S)={\rad}_C(S)=1.$$

\noindent
Indeed, for any  $x,y\in S$,  we have $\|x-y\|\leq \max\{1,\|x-y\|_{\ell_\infty(\N)}\}\leq 1$, so $\diam (S)\leq 1$. The vectors $(1,0,\dots)$ and $(0,1,0,\dots)$  are in $S$ and their distance from one another  equals $1$, so $\diam (S) =1$.
Now fix $y\in\cX$. If $\|y\|\geq 1$,  we have $\sup_{x\in S}\|y-x\|\geq \|y-0\|\geq 1$. 
On the other hand, if $\|y\|\leq 1$, then for any $\e>0$ there exists a coordinate $y_{j_0}$,  such that $|y_{j_0}|\leq \e$. Let $z\in S$ have $j_0$-th coordinate equal to $1$ and the rest of coordinates equal to $0$. Then, $ \|y-z\|\geq 1-\e$, so we get $\rad(S)=1$.  Then, from \eref{A1} we also have $ {\rad}_C(S)=1$.
\end{example}

\noindent
The following examples show that the same behavior can happen in classical spaces without modifying norms.

 \begin{example}
Let $\cX=c_0$ and $S:=\{x\ :\ x_j\geq 0 \mbox{ and } \sum_{j=1}^\infty x_j=1\}$.
Then, we have that 
$$\rad(S)={\rad}_C(S)=\diam (S)=1.
$$
 We can modify this example by taking  $\cX=\ell_p(\N)$, $1<p<\infty$,  and $S$ as above. In this case,  $\diam(S)=2^{1/p}$ and ${\rad}_C(S)=\rad(S)=1$.
\end{example}
\begin{example}
Let $\cX=L_1([0,1])$ and  $S:=\{f\in L_1([0,1]) :\ \int_0^1f=\int_0^1|f|=1\}$.
Then $\diam(S)=2={\rad}_C(S)$. However, by taking the ball centered at zero, we see that  $\rad(S)=1$. 
\end{example}

\begin{example}\label{Exnonsymmetric}
Let  $\cX:=\R^3$ with the   $\|\cdot \|_{\ell_\infty(\R^3)}$ norm. Let us consider the simplex
$$T:=\{x=(x_1,x_2,x_3)\ :\ \|x\|_\infty\leq 1 \mbox{  and  }x_1+x_2+x_3=2\}$$ 
with vertices $(1,1,0)$, $(1,0,1)$, and  $(0,1,1)$. We have 
$$
\diam(T)=1, \quad \rad(T)=\frac12, \quad {\rad}_C(T)=\frac23.
$$

\noindent
Indeed, since $T$ is the convex hull of its vertices, any  point in $T$ has coordinates in $[0,1]$, and hence the distance between
any two such points is at most one.   Since the vertices are at distance one from each other, we have that $\diam(T)=1$.
 It follows from  \eref{A1} that  $\rad(T)\ge 1/2$. Note that  the ball with center $(\frac12,\frac12,\frac12)$ and radius $1/2$ contains $T$,
  and so $\rad(T)=1/2$.  
 Given any point $z\in T$ which is a potential center of the restricted Chebyshev ball for $T$,
at least one of the coordinates of $z$ is at least $2/3$ (because $z_1+z_2+z_3=2$),
 and thus has distance at least $2/3$ from one of the vertices of $T$. 
  On the other hand, the ball with center $(\frac 2 3, \frac 2 3, \frac 2 3)\in T$ 
 and radius $\frac 2 3$ contains $T$.
\end{example}

\subsection{Connection to approximation sets and measurements}
\label{aprsets}

\noindent
The examples in this section are directed at showing that the  behavior, observed in \S\ref{exmp}, 
can occur even when the sets $S$ are described through measurements.
The next example is  a modification of Example \ref{Exnonsymmetric}, and  the set under 
consideration is of the form $\cK_w$,  
where $\cK$ is an approximation set.

\begin{example}
\label{example2}
We   take $\cX:=\R^4$ with  the $\|\cdot \|_{\ell_\infty(\R^4)}$ norm and define $V$ as  the one dimensional subspace 
spanned by $e_1:=(1,0,0,0)$. We consider the approximation set  
$$
\cK=\cK(1,V):=\{x\in\R^4: \dist(x,V)\le 1\},
$$
and   the measurement operator   $M(x_1,x_2,x_3,x_4)=(x_1,x_2+x_3+x_4)$.
Let us  now take the measurement  $w=(0,2)\in \R^2$ and look at $\cK_w$. 
 Since
$$
\cK= \{(t,x_2,x_3,x_4):\  t\in \R,\ \max_{2\le j\le 4}|x_j|\leq 1\}, \quad 
\cX_w= \{(0,x_2,x_3,x_4): \   x_2+x_3+x_4=2\},
$$
we infer that $\cK_w=\cX_w\cap\cK=\{ (0, x) \ : \  x\in \ T\}$,  where $T$ is the set from Example {\rm\ref{Exnonsymmetric}}. 
Thus,  we have 
$$\diam(\cK_w)=1, \ \ \ \ \ \ \  \rad(\cK_w)=\tfrac12,\ \ \ \ \ \ \ \    {\rad}_C(  \cK_w)=\tfrac23.$$
\end{example}

The following theorem shows that   any example for general sets $S$  can be transferred to the setting of interest to us, 
where the sets  are  of the form $\cK_w$ with $\cK$ being an approximation set.    

 \begin{theorem} 
 \label{theorem:asets}
 Suppose $X$ is a Banach space and $K\subset  X$ is a non-empty, closed and  convex subset of the closed unit ball $U$ of $X$.
Then, there exists a Banach space $\cX$, a finite dimensional subspace $V$, a measurement operator $M$, and a measurement $w$,
 such  that for the approximation set $\cK:=\cK(1,V)$, we have
$$
\diam(\cK_w)=\diam(K), \quad
\rad(\cK_w)=\rad(K), \quad
{\rad}_C(\cK_w)={\rad}_C(K).
$$
\end{theorem}

\noindent
{\bf Proof :}
Given $X$, we first define
$Z:=X\oplus \R:=\{(x,\alpha):\ x\in X, \ \alpha\in \R\}.$
Any norm on $Z$ is determined by describing its unit ball, which can be taken as 
any  closed,  bounded, centrally symmetric convex set.  We take  the set  $\Omega$  
to be  the convex hull of the set $(U,0)\cup (K,1)\cup(-K,-1)$. 
  Since $K\subset U$, it follows that a point of the form $(x,0)$ is in  $ \Omega$
  if and only if  $\|x\|_X\leq 1$.   Therefore,  for any $x\in X$,  
  \begin{equation}\label{123}
  \|(x,0)\|_Z=\|x\|_X.
   \end{equation}
Note also that  for any point   $ (x,\alpha)\in\Omega$,    we  have   $\max\{\|x\|_X,|\alpha|\}\leq 1$, and thus
   \begin{equation}\label{124}
   \max\{\|x\|_X,|\alpha|\}\leq \|(x,\alpha)\|_Z.
\end{equation}
It follows  from \eref{123} that for any  $x_1,x_2\in X$,  we have $\|(x_1,1)-(x_2,1)\|_Z=\|x_1-x_2\|_X$.  
 Now  we define
$\tilde K:=(K,1)\subset Z$.   Then, we have
$$\diam(\tilde K)_Z=\diam (K)_X, \quad
{\rad}_C(\tilde K)_Z={\rad}_C(K)_X.
$$
Clearly,  $\rad(\tilde K)_Z\leq \rad(K)_X$. On the other hand, for each $(x',\alpha)\in Z$,  we have
$$
\sup_{(x,1)\in \tilde K}\|(x,1)-(x',\alpha)\|_Z=\sup_{x\in K}\|(x-x',1-\alpha)\|_Z\geq \sup_{x\in K}\|x-x'\|_X
\geq \rad(K)_X,
$$
where the next to  last inequality uses \eref{124}.   Therefore, we   have $\rad(K)_X=\rad(\tilde K)_Z$.  
Next, we  consider the functional $\Phi\in Z^*$,  defined by  
 $$
  \Phi(x,\alpha)=\alpha.  
 $$
It follows from \eref{124} that  it  has norm one and 

\begin{equation} \label{Ex1}
\{z\in Z :   \Phi(z)=1,\ \|z\|_Z\le 1\}=\{(x,1)\in\Omega\} =\{(x,1): \ x\in K\}=\tilde K,
\end{equation}
where the next to the  last equality uses the fact that  a point of the form $(x,1)$ is in  $\Omega$ if and only if  $x\in K$.
We next define   the space  $\cX=Z \oplus \R:=\{(z,\beta):\ z\in Z, \beta\in\R\}$,  with the norm
$$
\|(z,\beta)\|_\cX:=\ \max\{\|z\|_Z,|\beta|\}.
$$
Consider the subspace $V=\{ (0,t) :\ t\in \R\}\subset  \cX$. 
If we take $\e=1$, then the approximation set $\cK=\cK(1,V)\subset \cX$  is
$\cK=\{(z,t):\ t\in\R, \ \|z\|_Z\le 1\}$.
We  now take the measurement operator $M(z,\beta)= (\beta,\Phi(z))\in\R^2$ and the  measurement $w=(0,1)$  which gives
$\cX_w=\{(z,0)\ :\  \Phi(z)=1\} $.
Then, because of \eref{Ex1}, we have
$$
\cK_w= \{(z,0): \Phi(z)=1,\ \|z\|_Z\le 1\}  = (\tilde K,0).
$$
   As above,  we prove that 
$$
\diam((\tilde K,0))_\cX=\diam (\tilde K)_Z, \quad {\rad}_C((\tilde K,0))_\cX={\rad}_C(\tilde K)_Z,  \quad \rad((\tilde K,0))_\cX=\rad(\tilde K)_Z,
$$
 which completes the proof of the theorem. \hfill $\blacksquare$

\section{A description of  algorithms via liftings} 
\label{sec:liftings}
In this section, we show that  algorithms for the optimal recovery problem can be described by what are called liftings in 
the theory of Banach spaces.   We place ourselves in the setting that $S$ is any subset of $\cX$, and we wish to recover 
the elements in $S_w$ for each measurement $w\in\R^m$.    That is, at this stage, we do not require that $S$ is an approximation set. 
Recall that given the measurement functionals $l_1,\dots,l_m$ in $\cX^*$,   the  linear operator $M:\cX\rightarrow \R^m$
is defined as
$$
M(f):=(l_1(f),\dots,l_m(f)),\quad f\in \cX.
$$
 Associated to $M$ we have the null space  
$$
 \cN:=\ker M=\{f\in \cX: \ M(f)=0\}\subset \cX,
$$
  and
$$
 \cX_w:= M^{-1}(w):=\{f\in\cX: M(f)=w\}.
$$
Therefore $\cX_0=\cN$.  Our goal is to recover the elements in     $S_w=\cX_w\cap \cS$.

\begin{remark}
\label{equivinfo}
Let us note that if in place of $l_1,\dots, l_m$, we use functionals $l_1',\dots,l_m'$ which span the same space  $L$ in $X^*$,
then the information about $f$ contained in $M(f)$ and $M'(f)$ is exactly the same,  and so the recovery problem is identical.
For this reason,  we can choose any spanning set of linearly independent functionals in defining $M$ and obtain exactly 
the same recovery problem.  Note that, since these  functionals are linearly independent,  $M$ is a linear mapping from 
$\cX$ onto $\R^m$.
\end{remark}
   
 We begin by analyzing the  measurement operator $M$.   We   introduce  the following norm on $\R^m$ induced by $M$
\be
\label{defnormM}
\|w\|_M=\inf_{x\in \cX_w}\| x\|,
\ee
and  consider the quotient space $\cX/\cN$.  Each element in $\cX/\cN$ is a coset    $\cX_w$, $w\in\R^m$.   The quotient norm on this space is
given by
\be
\label{quotientnorm}
\|\cX_w\|_{\cX/\cN}=\|w\|_M.
\ee
The mapping $M$ can be interpreted as mapping $\cX_w\rightarrow w$ and,  in view of \eref{quotientnorm},  
is an isometry from $\cX/\cN$ onto $\R^m$ under the norm $\|\cdot\|_M$.
\vskip .1in

\noindent
{\bf Lifting Operator:}  
{\it A  lifting  operator $\Delta$ is a mapping from $\R^m$ to $\cX$ which assigns to each $w\in\R^m$ an 
element from the coset $\cX_w$, i.e., a representer of the coset.   }
\vskip .1in
 
Recall that any algorithm $A$  is a mapping   from $\R^m$ into $\cX$.  We would like  the mapping $A$ for our recovery problem 
to send $w$ into an element of $S_w$, provided  $S_w\neq \emptyset$, since then we would know that $A$ is nearly optimal 
(see Theorem \ref{theorem:nearoptimal}) up to the constant $2$.  
 So, in going further, we consider only algorithms $A$ which take $w$ into $\cX_w$.   At this stage we are not yet invoking our desire that $A$
 actually maps into $S_w$, only that it maps into $\cX_w$.
\vskip .1in

\noindent
{\bf Admissible Algorithm:}
 {\it We say that an algorithm $A:\R^m \rightarrow \cX$  is admissible if, for each $w\in\R^m$, $A(w)\in \cX_w$. }
\vskip .1in

Our interest in lifting operators is because  any  admissible algorithm  $A$ is a lifting $\Delta$,  and the performance of 
such an  $A$ is related to the norm of $\Delta$.   
A natural lifting, and the one with minimal norm   $1$,  would be one which maps $w$ into an element of minimal 
norm in $\cX_w$.  Unfortunately, in general, no such minimal norm element exists, as the following illustrative example shows.

\begin{example}
\label{example3}
 We consider the space $\cX=\ell_1(\N)$ with the  $\|\cdot\|_{\ell_1(\N)}$ norm, and a collection of vectors 
$h_j\in \R^2$, $j=1,2,\dots$, with  $\|h_j\|_{\ell_2(\R^2)}=\langle h_j, h_j\rangle=1$, which are dense on the unit circle.  
We define the measurement operator $M$ as
$$
M(x):=\sum_{j=1}^\infty x_jh_j \in \R^2.
$$
 If follows from the definition of $M$ that for every $x$ such that $M(x)=w$, we have $\|w\|_{\ell_2(\R^2)}\le \|x\|_{\ell_1(\N)}$, and 
thus  $\|w\|_{\ell_2(\R^2)}\le \|w\|_M$.  In particular,  for every $i=1, 2, \ldots$, 
$$
1=\|h_i\|_{\ell_2(\R^2)}\le \|h_i\|_M\leq \|e_i\|_{\ell_1(\N)}=1,
$$
since $M(e_i)=h_i$, where $e_i$ is the $i$-th coordinate vector in $\ell_1(\N)$. So, we have that $\|h_i\|_{\ell_2(\R^2)}=\|h_i\|_M=1$. 
Since the $h_i$'s  are dense on the unit circle, every $w$ with Euclidean norm one satisfies $\|w\|_M=1$.
Next, we consider any  $w\in \R^2$,  such that $\|w\|_{\ell_2(\R^2)}=1$,  $w\neq h_j$, $j=1,2,\ldots$.  
If $w=M(x)=\sum_{j=1}^\infty x_jh_j$, then
$$
1=\langle w, w\rangle =\sum_{j=1}^\infty x_j\langle w,h_j\rangle.
$$
Since the $|\langle w,h_j\rangle|<1$, we must have  $\|x\|_{\ell_1(\N)}>1$.  Hence, $\|w\|_M$ is not assumed by any element
$x$  in the coset $\cX_w$.  This also shows there is no lifting $\Delta$ from $\R^2$ to $\cX/\cN$ with norm one.
\end{example}

While the above example shows that norm one liftings may not exist for a general Banach space $\cX$, there is a  
classical theorem of Bartle-Graves  which states  that there are  continuous liftings $\Delta$ with norm $\|\Delta\|$ 
as close to one as we wish (see \cite{BG, BL, RS}).  In our setting,  this theorem can be stated as follows.

\begin{theorem}[Bartle-Graves]\label{BartleGraves} Let $M:\cX\rightarrow \R^m$ be a measurement operator.  
For every $\eta>0$, there exists a   map $\Delta:\R^m\rightarrow \cX$, such that
\begin{itemize}
\item $\Delta$ is continuous.
\item $\Delta(w)\in \cX_w,\quad w\in\R^m$.
\item for every  $\lambda >0$,  we have  $\Delta(\lambda w)=\lambda\Delta(w)$.
\item $\|\Delta(w)\|_\cX\leq (1+\eta)\|w\|_M,\quad w\in \R^m$.
\end{itemize}
 \end{theorem}

Liftings are closely related to projections.  If $\Delta $ is a linear lifting,  then its range $Y$ is a subspace  of $\cX$ of dimension $m$,  
for which we have the following.
\begin{remark}   For any fixed constant $C$, there is  a linear  lifting $\Delta:\R^m\rightarrow \cX$ with norm $\le C$   if and only if there exists a linear projector $P$ from  $\cX$ onto a subspace $Y\subset \cX$ with  $\ker (P)=\cN$ and  $\|P\|\leq C$. 
\end{remark}
\noindent
{\bf Proof:}
  Indeed, if $\Delta$ is such a lifting then its range $Y$ is a finite dimensional subspace and $P(x):=\Delta(M(x))$
defines a projection  from $\cX$ onto $Y$ with the mentioned properties.  On the other hand, given such a $P$ and $Y$,   notice that any two elements in $M^{-1}(w)$ have the same image under $P$, since  the kernel of $P$ is $\cN$.   Therefore, we can define the lifting $\Delta(w):=P(M^{-1}(w))$, $w\in\R^m$, which has norm at most $C$. 
  \hfill $\blacksquare$

The above results are for an arbitrary Banach space.  If we put more structure on $\cX$, then we can guarantee the existence of
a continuous lifting with norm one  (see  \cite[Lemma 2.2.5]{BL}). 
\begin{theorem}
\label{theorem:lifting}
If the Banach space $\cX$ is uniformly convex,  then for each $w\in\R^m$, there is a unique
$x(w)\in\cX_w$,  such that
\be
\label{tl1}
\|x(w)\|_\cX= \inf_{x\in\cX_w}\|x\|=:\|w\|_M.
\ee
The mapping $\Delta:w\rightarrow x(w)$ is a continuous  lifting of norm one.
\end{theorem}

\noindent
{\bf Proof:}  Fix $w\in\R^m$ and let $x_j\in \cX_w$, $j\ge 1$, be such that $\|x_j\|\to \|w\|_M$.  Since $\cX$ is uniformly convex, by  weak  compactness, there is
a subsequence of $\{x_j\}$ which, without loss of generality, we can take as $\{x_j\}$ such that $x_j\to x \in \cX$ weakly.  
It follows that $\lim_{j\to\infty}l(x_j)=l(x)$ for all $l\in \cX^*$.  Hence $M(x)=w$, and therefore $x\in \cX_w$.   Also, if 
$l$ is a norming functional for $x$, i.e. $\|l\|_{\cX^*}=1$ and  $l(x)=\|x\|$,  then
$$
\|x\|=l(x)=\lim_{j\to\infty}l(x_j)\le  \lim_{j\to\infty}\|x_j\|=\|w\|_M,
$$
 which shows the existence in \eref{tl1}.  
To see that  $x=x(w)$  is unique,  we assume $x'\in \cX_w$ is
another element with $\|x'\|=\|w\|_M$.  Then $z:=\frac{1}{2}(x+x')\in\cX_w$,  and by uniform  convexity $\|z\|<\|w\|_M$,  which is an obvious contradiction.  This shows that there is an $x=x(w)$ satisfying \eref{tl1}, and it is unique.

\noindent
  To see that $\Delta$ is continuous, let
$w_j\to w$ in $\R^m$ and let $x_j:=\Delta(w_j)$ and $x:=\Delta(w)$.    Since we also have that 
$\|w_j\|_M\rightarrow \|w\|_M$, it follows from  the minimality of $\Delta(w_j)$ that $\|x_j\|\to \|x\|$.
 If $w=0$, we have  $x=0$, and thus we have convergence in norm. In what follows, we assume that 
$w\neq 0$.    Using weak compactness (passing to a subsequence),  we can assume that $x_j$ converges weakly to some  $\bar x$. So, we have  $w_j=M(x_j)\rightarrow M(\bar x)$,  which gives $M(\bar x)=w$. Let $\bar l\in \cX^*$ be a norming functional for $\bar x$.  Then, we have that 
$$
\|\bar x\|=\bar l(\bar x)=\lim_{j\to\infty} \bar l(x_j)\le \lim_{j\to\infty} \|x_j\|=\|x\|,
$$
and therefore  $\bar x=x$ because of the definition of $\Delta$. We want to show that $x_j\to x$ in norm. If this is not the case, we can find a subsequence, which we again denote by $\{x_j\}$,  such that $\|x_j-x\|\ge \e>0$, $j=1,2,\dots$, for some $\e>0$.     It follows 
from the  uniform convexity  that $\|\frac{1}{2}(x_j+x)\|\le \max\{\|x_j\|,\|x\|\}\alpha$ for all $j$, with $\alpha<1$ a fixed constant.  Now, let $l\in \cX^*$ be a norm one 
functional,  such that $l(x)=\|x\|$. Then,  we have
$$
2\|x\|= 2l(x)=\lim_{j\to\infty} l(x_j+x) \le \lim_{j\to\infty} \|x_j+x\|\le2 \|x\|\alpha ,
$$
which gives $\alpha\geq 1$ and is the desired contradiction.\hfill $\blacksquare$

\begin{remark}
\label{remark:Brown}
The paper {\rm\cite{Br}} gives an example of a strictly convex, reflexive Banach space $\cX$ and a measurement map $M:\cX\rightarrow \R^2$,
   for which there is no continuous norm one lifting $\Delta$.  Therefore, the above theorem would not hold under the slightly milder assumptions on $\cX$ being strictly convex and reflexive (in place of uniform convexity).
\end{remark}


\section{A priori estimates for the radius of $\cK_w $}
\label{sec:radius}
In this section, we discuss estimates for  the radius of $\cK_w $ when $\cK=\cK(\e,V)$ is an approximation set.  The main result we shall obtain is that the global optimal recovery error $R(\cK)$
is determined a priori (up to a constant factor $2$) by the angle between the null space $\cN$ of the measurement map $M$
and the approximating space $V$ (see (iii) of Theorem \ref{theorem:main} below).

\begin{remark}
\label{radiusremark} 
Note the following simple observation.  

\noindent
{\rm (i)}  If    $  \cN \cap V\neq \{0\}$,    then for any $0\neq \eta \in \cN\cap V$,    and any  $x\in \cK_w$,  the  line $x+t \eta$, $t\in\R$, is contained in $\cK_w$, and therefore there is no finite ball $B(a,r)$ which contains $\cK_w$.  Hence $\rad(\cK_w)=\infty$.
 
 \noindent
{\rm (ii)} If  $\cN\cap V=\{0\}$, then $n=\dim V \leq \codim \cN=\rank M=m$ and therefore $n\leq m$.  In this case $\rad(\cK_w)$ is finite for all $w\in\R^m$.
\end{remark}

\noindent
{\bf Standing Assumption:} {\it In view of this remark, the only interesting case is {\rm (ii)}, and therefore we assume that 
$\cN\cap V=\{0\}$  for the remainder of this paper.}   
\vskip .1in

For (arbitrary) subspaces $X$ and $Y$ of a given Banach space $\cX$, we recall the   angle $\Theta$ between $X$ and $Y$, 
defined as
$$
 \Theta(X,Y):= \inf_{x\in X} \frac{ \dist(x,Y)}{\|x\|}.
$$
We are more interested in $\Theta(X,Y)^{-1}$,  and so accordingly, we  define
\be
\label{defmu}
\mu(X,Y):= \Theta(X,Y)^{-1}=\sup_{x\in X}\frac{\|x\|}{\dist(x,Y)}=\sup_{x\in X,y\in Y}\frac{\|x\|}{\|x-y\|}.
\ee
Notice that $\mu(X,Y)\geq 1$.

 \begin{remark}
Since $V$ is a finite dimensional space and $\cN\cap V=\{0\}$,  we have $\Theta(\cN,V)>0$.  Indeed, otherwise there exists a 
sequence $\{\eta_k\}_{k\ge 1}$ from $\cN$ with $\|\eta_k\|=1$ and a sequence  $\{v_k\}_{k\ge 1}$ from $ V$, 
such that $\|\eta_k-v_k\|\to 0$, $k\to\infty$. We can assume $v_k$ converges to $v_\infty$,  but then also $\eta_k$ 
converges to $v_\infty$,  so $v_\infty \in \cN\cap V$ and $\|v_\infty\|=1$, which is the desired contradiction to $\cN\cap V=\{0\}$.
\end{remark}

Note that,   in general,  $\mu$ is not symmetric, i.e.,  $\mu(Y,X)\neq \mu(X,Y)$.
However, we do have the following comparison.

 \begin{lemma}
\label{mu}
For arbitrary subspaces $X$ and $Y$ of a given Banach space $\cX$,  such that $X\cap Y=\{0\}$, we have  
\be
\label{ineq}
\mu(X,Y)\leq 1+\mu(Y,X)\le 2 \mu(Y,X).
\ee
\end{lemma}
\noindent
{\bf Proof:} For each $x\in X$ and $y\in Y$ with $x\neq 0$,  $x\neq y$, we have 
$$
\frac{\|x\|}{\|x-y\|}\le  \frac{\|x-y\|+\|y\|}{\|x-y\|}=1+  \frac{\|y\|}{\|x-y\|}\le 
 1+ \mu(Y,X).
$$
Taking a supremum over $x\in X,y\in Y$, we arrive at  the first inequality in \eref{ineq}.  The second inequality follows because $\mu(Y,X)\ge 1$.  $\hfill \blacksquare$

The following  lemma records some properties of $\mu$ for our setting in which $Y=V$ and $X=\cN$ is the null space of $M$.

\begin{lemma}
\label{comp}
Let  $\cX$ be any Banach space, $V$ be any finite dimensional subspace  of $ \cX$ with $\dim(V)\le m$,   and  $M:\cX\to \R^m$ be any measurement operator.   Then,  for the null space $\cN$ of $M$,  we have the following.

\noindent
{\rm (i)} $
\mu(V,\cN)=\|M_V^{-1}\|,
$

\noindent {\rm (ii)} $
\mu(\cN,V) \leq 1+\mu(V,\cN)=1+\|M_V^{-1}\|\le 2\|M_V^{-1}\|,
$

\noindent
 where $M_V$ is the restriction of the measurement operator $M$ on $V$ and $M_V^{-1}$ is its inverse.
\end{lemma}
{\bf Proof:}
The statement (ii) follows from (i) and Lemma \ref{mu}.  To prove (i), we see from the definition of $\|\cdot\|_M$ given in \eref{defnormM}, we have
 \begin{eqnarray}
 \nonumber
\|M_V^{-1}\|= \sup_{v\in V}
\frac{\|v\|}{\|M_V(v)\|_M} =\sup_{v\in V}
\frac{\|v\|}{\dist(v,\cN)}=\mu(V,\cN),
\end{eqnarray}
as desired. \hfill $\blacksquare$

We have the following simple,  but important theorem.
\begin{theorem}
\label{theorem:main}
Let  $\cX$ be any Banach space, $V$ be any finite dimensional subspace  of $ \cX$,  $\e>0$, 
 and  $M:\cX\to \R^m$ be any measurement operator.   Then, for the set $\cK=\cK(\e,V)$, we have the following

\noindent
 {\rm (i)} For any $w\in \R^m$, such that $w=M(v)$ with $v\in V$, we have   
$$
\rad (\cK_w)=\e\mu(\cN,V).
$$
\noindent
{\rm (ii)} For any $w\in \R^m$, we have
$$
\rad (\cK_w)\le  2\e\mu(\cN,V).
$$
\noindent
{\rm (iii)} We have
$$
\e\mu(\cN,V)\le R(\cK)\le 2\e\mu(\cN,V).
$$
\end{theorem}
 
\noindent
 {\bf Proof:} First,  note that  $\cK_0=\cK\cap \cN$  is centrally symmetric and convex and likewise $\cK_0(\e,V)$ is also centrally symmetric and convex.   Hence, from  Proposition \ref{prop:centrallysymmetric}, we have that the smallest ball containing 
 this set is centered at $0$ and   has radius
 \be 
 \label{main}
\rad(\cK_0(\e,V))=   \sup\{ \|z\| :\ z\in \cN,  \dist(z,V)\leq \e\}=
\e\mu(\cN,V).
\ee
      Suppose now that  $w=M(v)$ with $v\in V$.  Any  $x\in\cK_w$ can be written as $x=v+\eta$ with $\eta\in\cN$ if and only if $\dist(\eta,V)\le \e$.   Hence $\cK_w=v+\cK_0$ and (i) follows. 
      
       For the proof of  (ii), let $x_0$ be any point in $\cK_w$.  Then, any other $x\in\cK_w$ can be written as $x=x_0+\eta$.  Since $\dist(x,V)\le \e$, we have $\dist(\eta,V)\le 2\e$.  Hence 
 $$
 \cK_w \subset\, x_0+\cK_0(2\e,V),
 $$
  which from  \eref{main} has radius  $2\e \mu(\cN,V)$. Therefore,  we have proven  (ii).  Statement (iii) follows from the definition 
  of $R(\cK)$ given in \eref{defR}.  
  \hfill $\blacksquare$

 Let us make a few comments about Theorem \ref{theorem:main} viz a viz the results in \cite{BCDDPW1}  
 (see Theorem 2.8 and Remark 2.15 of that paper) for the case  when $\cX$ is a Hilbert space.  In the latter case, it was shown 
 in \cite{BCDDPW1}  that the same result  as (i) holds,  but in the case of (ii), an exact computation of $\rad(K_w)$ was given with the constant $2$ replaced by a number (depending on $w$) which is less than one.  
  It is probably impossible to have an exact formula for $\rad(K_w)$   in the case of a general Banach space. However,  we show in the appendix that when  $\cX$ is uniformly convex and uniformly smooth,  we can improve on the constant appearing in   (ii) of Theorem \ref{theorem:main}.

\section{ Near optimal  algorithms  }
\label{sec:algorithms}
In this section, we discuss the form of admissible algorithms for optimal recovery,   and expose what properties
these algorithms need in order to be optimal   or near  optimal on the classes $\cK_w$ when $\cK=\cK(\e,V)$.    
Recall that any algorithm$A$ is a mapping $A:\R^m\rightarrow \cX$. Our goal is  to have  $A(w)\in \cK_w$ 
for each $w\in \R^m$,  for which
$\cK_w\neq \emptyset$, since, by Theorem \ref{theorem:nearoptimal}, this would guarantee that the algorithmic error
\begin{eqnarray}
\label{plo1}
\sup_{x\in \cK_w}\|x-A(M(x))\|\le   2\rad(\cK_w),
\end{eqnarray}
and hence up to the factor $2$ is optimal. 
In this section, we shall not be concerned about computational issues that arise in  the numerical implementation 
of the algorithms we put forward.   Numerical implementation issues  will be discussed in the section that follows.

Recall that by $M_V$ we denoted  the restriction of $M$ to the space $V$.   By our {\bf Standing Assumption},  
$M_V$ is invertible,  and hence $Z:=M(V)=M_V(V)$ is an $n$-dimensional subspace of $\R^m$.   Given $w\in \R^m$,  we consider
its error of best approximation from $Z$ in $\|\cdot\|_M$,  defined by
$$
E(w):=\inf_{z\in Z}\|w-z\|_M.
$$
Notice that whenever $w=M(x)$,  from the definition of the norm $\|\cdot\|_M$, we have
\be
\label{equiv}
E(w)=\dist(w,Z)_M=\dist(x,V\oplus \cN)_\cX\le \dist(x,V)_\cX.
\ee

While there is always a best approximation $z^*=z^*(w)\in Z$ to $w$,  and it is unique when the norm is strictly convex, 
for a possible ease of  numerical implementation,
we consider other non-best  approximation maps.   
We say a  mapping $\Lambda:\R^m\mapsto Z$ is {\it near best} with constant $\lambda\geq 1$,   if
\be
\label{nbZ}
\|w-\Lambda(w)\|_M\le \lambda E(w),\quad w\in\R^m.
\ee
Of course, if $\lambda=1$, then $\Lambda$ maps $w$ into a best approximation of $w$ from $Z$.

Now, given any lifting $\Delta$ and any near best approximation map $\Lambda$, 
we consider the mapping
\be
\label{algmap}
A(w):= M_V^{-1}(\Lambda(w))+\Delta(w-\Lambda(w)),\quad w\in \R^m.
\ee
Clearly,  $A$ 
 maps $\R^m$ into $\cX$,  so that it is an algorithm.  It also has the property that $A(w)\in \cX_w$, 
 which means that it is an admissible algorithm.   Finally,   by  our construction, whenever $w=M(v)$ 
 for some $v\in V$, then    $\Lambda(w)=w$, and so  $A(w)=v$.   
  Let us note some  important properties of such an  algorithm   $A$.  
 \begin{theorem}
 \label{theorem:ba}  
 Let  $\cX$ be a Banach space, $V$ be any finite dimensional subspace  of $ \cX$,  $\e>0$,     
 $M:\cX\to \R^m$ be any measurement operator  with a null space $\cN$, and 
 $\cK=\cK(\e,V)$ be an approximation set.   Then, for any lifting $\Delta$ and any  near best approximation 
 map $\Lambda$ 
 with constant $\lambda\geq 1$, the algorithm $A$, defined in {\rm\eref{algmap}},  has the following properties
 
 \noindent
{\rm (i)}\,\,$A(w)\in\cX_w,\quad w\in \R^m$.
 
 \noindent
{\rm (ii)}\,\,$\dist(A(M(x)),V)\le \lambda\|\Delta\| \dist(x,V)_\cX,\quad x\in\cX.$
 
  \noindent
{\rm (iii)}\,\, if $\|\Delta\|=1$ and $\lambda=1$, then  $A(M(x))\in\cK_w$, whenever $x\in\cK_w$.
  
 \noindent
{\rm (iv)}\,\, if $\|\Delta\|=1$ and $\lambda=1$, then the algorithm $A$ is near optimal with constant $2$, i.e. for any $w\in\R^m$,
 \be
 \label{tnb}
 \sup_{x\in \cK_w}\|x-A(M(x))\|\le 2 \rad(\cK_w).
  \ee
 
 \noindent
{\rm (v)} \,\, if $\|\Delta\|=1$ and $\lambda=1$, then the algorithm $A$ is also near  optimal  in the sense of minimizing $R(\cK)$,  and
$$
\sup_{w\in\R^m}\sup_{x\in \cK_w} \|x-A(M(x))\|\le 2R(\cK).
$$

 \noindent
{\rm (vi)}\,\, if $\|\Delta\|=1$ and $\lambda=1$, then the algorithm  $A$ has the a priori performance bound 
\be
\label{algper}
\sup_{w\in\R^m}\sup_{x\in \cK_w} \|x-A(M(x))\|\le 4\e\mu(\cN,V).
\ee

 \end{theorem}
 
\noindent
{\bf Proof:}   We have already noted that (i) holds.   To prove (ii), let $x$ be any element in $\cX$.   
Then $M_V^{-1}(\Lambda(M(x)))\in V$, and therefore
$$
\dist(A(M(x)),V)_\cX\le \|\Delta \| \|M(x)-\Lambda(M(x))\|_M\le    \|\Delta \|\lambda E(M(x))\le \|\Delta\|\lambda \dist(x,V),
$$
where the  first inequality uses \eref{algmap}, the second inequality uses \eref{nbZ},  and 
the last equality uses \eref{equiv}.    The statment  (iii) follows from (i) and (ii), since whenever $x\in \cK_w$, then $\dist(x,V)_\cX\le \e$. 
The statement  (iv) follows from (iii) because of Theorem \ref{theorem:nearoptimal}.
The estimate  (v) follows from (iv) and the definition \eref{defR} of $R(\cK)$.  
Finally, (vi) follows from (v) and the a priori estimates of Theorem \ref{theorem:main}.  \hfill $\blacksquare$

 \subsection{Near best algorithms}
 \label{ssnearbest}
 In view of the last  theorem, from a theoretical point of view, the best choice for $A$ is to choose $\Delta$ with $\|\Delta\|=1$ and  $\Lambda$
 with constant $\lambda=1$.    When $\cX$ is uniformly convex, we  can always accomplish this theoretically, but there
  may be  issues in the numerical implementation.   If $\cX$ is a general Banach space, we can choose $\lambda=1$ 
  and $ \|\Delta\|$ arbitrarily close to one,  but  as in the latter case, problems in the numerical implementation may also arise.   
 In the next section, we discuss some of the numerical considerations in implementing an algorithm $A$ of the  form \eref{algmap}.
  In the case that $\lambda \|\Delta\|>1$,  we only know that
  $$
  \dist(A(M(x))),V)\le \lambda\|\Delta\| \e, \quad x\in \cK.
  $$
    It follows that $A(w)\in \cK_w(\lambda\|\Delta\|\e,V)$.   Hence, from \eref{plo1} and Theorem \ref{theorem:main}, we know that
  $$
  \sup_{x\in \cK_w}\|x-A(M(x)\|\le 4\lambda\|\Delta\|\e\mu(\cN,V).
 $$
This is only slightly worse than the a priori bound $4\e\mu(\cN,V)$ which we obtain when we know that $A(w)$ is in $\cK_w(\e,V)$.  In this case, 
the algorithm $A$ is near best for $R(\cK)$ with the constant  $4\lambda \|\Delta\|$.

\section{Numerical issues in implementing the algorithms $A$}
\label{sec:numerical}  
In this section, we address the main numerical issues in implementing algorithms of the form \eref{algmap}.  These are
\begin{itemize}
\item  How to compute $\|\cdot\|_M$ on $\R^m$?
\item  How to numerically construct near best approximation maps  $\Lambda$ for approximating 
the elements in $\R^m$ by the elements of $Z=M(V)$ in the norm
$\|\cdot\|_M$?
\item 
How to numerically construct lifting operators $\Delta$ with a controllable norm $\|\Delta\|$?
\end{itemize}
Of course, the resolution of each of these issues depends very much on the Banach space $\cX$, the subspace $V$, and the measurement functionals $l_j$, $j=1,\dots,m$.  In this section, we will consider general
principles   and see how these principles are implemented in three examples.
\vskip .1in

\noindent{\bf Example 1:}  $\cX=C(D)$,  where $D$ is a domain in $\R^d$, $V$ is any $n$ dimensional subspace of
$\cX$, and $M=(l_1,\ldots,l_m)$ consists of 
$m$ point evaluation functionals at distinct points $P_1,\dots,P_m\in D$, i.e., $M(f)=(f(P_1),\dots, f(P_m))$.

\vskip .1in
 \noindent{\bf Example 2:}  $\cX=L_p(D)$, $1\le p\le \infty$, where $D$ is a domain in $\R^d$, $V$ is any $n$ dimensional subspace of
$\cX$ and $M$ consists of the $m$   functionals
$$
l_j(f):=\int_D f(x)g_j(x)\,dx, \quad j=1,\dots,m,
$$
where
the functions $g_j$ have disjoint supports, $g_j\in L_{p'}$, $p'=\frac{p}{p-1}$,  and $\|g_j\|_{L_{p'}}=1$.  

\vskip .1in
\noindent{\bf Example 3:}  $\cX=L_1([0,1])$,  $V$ is any $n$ dimensional subspace of
$\cX$, and $M$ consists of the $m$   functionals
$$
l_j(f):=\int_0^1 f(t)r_j(t)\,dt, \quad j=1,\dots,m,
$$
where
the functions $r_j$ are the Rademacher functions
\be
\label{rademacher}
r_j(t):=  {\rm sgn } ( \sin 2^{j+1} \pi t), \quad  t\in [0,1], \quad j\ge 0.
\ee
The functions $r_j$  oscillate and have full support.  This example is not so important in typical data fitting scenarios, but it is important theoretically since, as we shall see, it  
has interesting features with regard to liftings.

\subsection{Computing $\|\cdot\|_M$.}
\label{ss:compM}
We assume that the measurement functionals $l_j$, $j=1,\dots,m$,  are given explicitly and  are linearly independent. Let  $L:=\span (l_j)_{j=1}^m\subset \cX^*$. Our strategy is to first compute the dual norm $\|\cdot\|_M^*$ on $\R^m$ by using the fact that the functionals $l_j$
are available to us.    Let $\alpha=(\alpha_1,\dots,\alpha_m)\in \R^m$ and consider its action as a linear functional.
 We have that 
 \begin{equation}\label{Quotdual2}
\|\alpha \|_M^*=\sup_{\|w\|_M=1}|\sum_{j=1}^m\alpha_jw_j|=\sup_{\|x\|_\cX= 1} \left|\sum_{j=1}^m\alpha_j l_j(x)\right|=\left\|\sum_{j=1}^m\alpha_j l_j\right \|_{\cX^*},
\end{equation}
 where we have used that if $\|w\|_M=1$, there exists $x\in\cX$, such that $M(x)=w$ and its norm is
arbitrarily close to one.
Therefore, we can express $\|\cdot\|_M$ as 
\begin{equation}\label{Quotdual1}
\|w\|_M=\sup\{|\sum_{j=1}^m w_j\alpha_j|\ \ :\ \  \|(\alpha_j)\|_M^*\leq 1\}.
\end{equation}

We consider these norms   to be computable since the space $\cX$ and the functionals $l_j$ are known to us.  Let us illustrate this
in our three examples.   In Example 1, for any $\alpha\in \R^m$, we have
$$
\|\alpha\|_M^*=\sum_{j=1}^m|\alpha_j|=\|\alpha\|_{\ell_1(\R^m)}, \quad \|w\|_M=\max_{1\le j\le m}|w_j|=\|w\|_{\ell_\infty(\R^m)}.
$$
 In  Example 2, we have $\cX=L_p$, and
$$
\|\alpha\|_M^*= \|\alpha\|_{\ell_{p'}(\R^m)}, \quad \|w\|_M= \|w\|_{\ell_p(\R^m)}.
$$
 In Example 3,  we have $\cX^*=L_\infty([0,1])$,  and from \eref{Quotdual2} we infer that
$\|\alpha\|_M^*=\|\sum_{j=1}^m \alpha_j r_j\|_{L_\infty([0,1])}$.   From the definition \eref{rademacher},  we see that the sum $\sum_{j=1}^k \alpha_j r_j$ is constant on each interval of the form $(s2^{-k}, (s+1)2^{-k})$ when $s$ is an integer.  On the other hand,  on such an interval, $r_{k+1}$ takes on    both of the  values $1$ and $-1$.  Therefore, by induction on $k$,  we get
$\|\alpha\|_M^*= \sum_{j=1}^m|\alpha_j|$.  Hence,
we have
$$
\|\alpha\|_M^*=\sum_{j=1}^m|\alpha_j|=\|\alpha\|_{\ell_1(\R^m)}, \quad\|w\|_M=\max_{1\le j\le m} |w_j|=
  \|w\|_{\ell_\infty(\R^m)}.
  $$

\subsection{Approximation maps}
\label{ss:approxoperators}
 
 Once  the norm $\|\cdot\|_M$ is numerically computable, the problem of finding a best or near best approximation map $\Lambda(w)$ to $w$
 in this norm becomes a standard problem in convex minimization.  For instance, in the examples from the previous subsection,
 the minimization is done in $\|\cdot\|_{\ell_p(\R^m)}$.
 Of course, in general, the performance of algorithms for such minimization depend on the properties of the unit ball of $\|\cdot\|_M$.  This ball is always convex, but in some cases it is uniformly convex and this  leads to faster convergence of the  iterative minimization algorithms
  and guarantees a unique minimum.

\subsection{Numerical liftings}
\label{ss:numliftings}

Given a prescribed null space $\cN$, a  standard way to find linear liftings from $\R^m$ to $\cX$ is to find a linear  
projection $P_Y$ 
from $\cX$ to a subspace $Y\subset \cX$  of dimension $m$ which 
has $\cN$ as its kernel.   We can find all  $Y$ that can be used in this fashion as follows.   We take  elements  $\psi_1,\dots,\psi_m$ from $\cX$, such that
$$
l_i(\psi_j)=\delta_{i,j},\quad 1\le i,j\le m,
$$
where $\delta_{i,j} $ is the usual Kronecker symbol.  In other words, $\psi_j$, $j=1,\dots,m$, is a dual basis to $l_1,\dots,l_m$.  
Then,  for $Y:=\span\{\psi_1,\dots,\psi_m\}$,  the projection
$$
P_Y(x)=\sum_{j=1}^ml_j(x)\psi_j, \quad x\in\cX,
$$
has kernel $\cN$.   We get a lifting corresponding to $P_Y$ by defining
\be
\label{liftingY}
\Delta(w):=\Delta_Y(w):=\sum_{j=1}^mw_j\psi_j.
\ee
This lifting is linear and hence continuous.   The important issue for us is its norm.  We see that
$$
\|\Delta\|=\sup_{\|w\|_M=1}\|\Delta(w)\|_\cX =\sup_{\|w\|_M=1}\|\sum_{j=1}^mw_j\psi_j\|_\cX= \sup_{\|x\|_\cX=1}\|\sum_{j=1}^ml_j(x)
\psi_j\|_\cX=\|P_Y\|.
$$
Here, we have used the fact that if $\|w\|_M=1$,  then there is an $x\in \cX$ with norm as close to one as we wish with $M(x)=w$.

It follows from the Kadec-Snobar theorem that we can always choose a $Y$ such that $\|P_Y\|\le \sqrt{m}$.  
In general, the $\sqrt{m}$ cannot be replaced by a smaller power of $m$.   However,  if $\cX=L_p$,  
then $\sqrt{m}$ can be replaced by $m^{|1/2-1/p|}$.  We refer the reader to Chapter III.B of \cite{Wbook} for a discussion of these facts.

In many settings, the situation is more favorable.   In the case of Example 1,   we can take for $Y$ 
the span of any norm one   functions $\psi_j$, $ j=1,\dots, m$,  such that $l_i(\psi_j)=\delta_{i,j}$, $1\le i,j\le m$.   
We can always take the $\psi_j$ to have disjoint supports,  and thereby get that $\|P_Y\|=1$.   
Thus, we get a linear lifting $\Delta$  with $\|\Delta\|=1$  (see \eref{liftingY}).    
This same discussion also applies to Example 2.

Example 3 is far more illustrative.   Let us first consider linear liftings $\Delta:\R^m\rightarrow L_1([0,1])$.   
It is well known (see e.g. \cite[III.A, III.B]{Wbook}) that we must have $\| \Delta\|\geq c\sqrt{m}$. 
A   well known, self-contained argument to prove  this is the following.     
Let  $e_j$, $j=1,\dots,m$, be the usual coordinate vectors in $\R^m$.   Then, the function  $ \Delta(e_j)=:f_j\in  L_1([0,1])$ and   
$\|f_j\|_{L_1([0,1])}\geq  \|e_j\|_M=\|e_j\|_{\ell_\infty(\R^m)}=1$.
Next, we fix   $\eta\in [0,1]$, and  consider for each fixed $\eta$
$$ 
\Delta(( r_1(\eta), \ldots,r_m(\eta)))=\Delta\left (\sum_{j=1}^mr_j(\eta)e_j\right)=\sum_{j=1}^m r_j(\eta) f_j(t).
$$
Clearly, $\|\Delta\|\ge \|\sum_{j=1}^m r_j(\eta) f_j\|_{L_1([0,1])}$ for each $\eta\in [0,1]$.
Therefore, integrating this inequality over $ [0,1]$ and using Khintchine's inequality with the best constant (see  \cite{Szarek}), we find
\begin{eqnarray*}
\|\Delta\|&\geq& \int_0^1 \|\sum_{j=1}^m r_j(\eta) f_j\|_{L_1([0,1])}\, d\eta =\int_0^1 \int_0^1|\sum_{j=1}^m r_j(\eta) f_j(t)|\,d\eta\, dt\\
&\geq&\frac1{\sqrt 2}\int_0^1\left(\sum_{j=1}^m f_j(t)^2\right)^{1/2} \, dt \geq \frac1{\sqrt 2} \int_0^1\frac1{\sqrt m}\sum_{j=1}^m|f_j(t)|\, dt\\
&=&\frac1{\sqrt 2} \frac1{\sqrt m}\sum_{j=1}^m \|f_j\|_{L_1([0,1])}\geq \frac{\sqrt m}{\sqrt 2},
\end{eqnarray*}
where the next to last inequality uses the Cauchy-Schwarz inequality.

 Even though linear liftings in Example 3 can never have a norm smaller than $
 \sqrt{m/2}$, we can construct nonlinear liftings which have norm one.    
 To see this, we  define such a lifting for   any $w\in\R^m$ with  $\|w\|_M=\max_{1\le j\le m}|w_j|=1$,  
 using the classical Riesz product construction.  Namely, for such $w$, we define 
\be
\label{Riesz}
\Delta(w):=\prod_{j=1}^m (1+w_jr_j(t))=\sum_{A\subset\{1,\dots,m\}} \prod_{j\in A}w_j r_j(t),
\ee
where we use the convention that    $ \prod_{j\in A}w_j r_j(t)=1$  when $A=\emptyset$.
Note that if $A\neq \emptyset $, then 
\be\label{126}
\int_0^1  \prod_{j\in A} r_j(t)\, dt=0.
\ee
Therefore,  $\int_0^1\Delta (w)\, dt =1=\|\Delta(w)\|_{L_1([0,1])}$,  because $\Delta(w)$ is a nonnegative function. 
To check that $M(\Delta(w))=w$, we first  observe  that 
\be
\label{integrals}
\left( \prod_{j\in A} r_j(t)\right) r_k(t)=\begin{cases}\prod_{j\in A\cup\{k\}}   r_j,& \mbox{ when }k\notin A,\\\\
\prod_{j\in A\setminus\{k\}}  r_j, & \mbox{ when }k\in A.
\end{cases}
\ee
Hence, from \eref{126} we see that  the only $A$ for which the integral of  the left hand side of \eref{integrals} is 
nonzero is when $A=\{k\}$.  This observation,  together with \eref{Riesz}  gives
$$
l_k(\Delta(w))=\int _0^1\Delta(w)r_k(t)\,dt = w_k, \quad 1\le k\le m,
$$
and therefore $M(\Delta(w))=w$.   We now define $\Delta(w)$ when $\|w\|_{\ell_\infty(\R^m)} \neq 1$ by
 \begin{equation}\label{Riesz1}
 \Delta(w)=\|w\|_{\ell_\infty}\Delta(w/\|w\|_{\ell_\infty}), \quad \Delta(0)=0.
\end{equation}
We have therefore proved that $\Delta$ is a lifting of norm one.


\section {Performance estimates for the examples}
\label{sec:performance}
 In this section, we consider the examples from \S\ref{sec:numerical}.  In particular, we determine
$\mu(\cN,V)$, which  allows us to give the global performance error for near  optimal 
algorithms for these examples.     We begin with the optimal algorithms in a 
Hilbert space, which is not one of our three examples, but is easy to describe.

\subsection{The case  when $\cX$ is a Hilbert space $\cH$}
\label{ss:exHilbert}
   This case was completely analyzed in   \cite{BCDDPW1}.
We summarize the results of that paper here in order to   point out that our algorithm is a direct extension 
of the Hilbert space case to the Banach space situation, and to compare this case with our examples 
in which  $\cX$ is not a Hilbert space. 
In the case $\cX$ is a Hilbert space,  the measurement functionals $l_j$  have the representation $l_j(f)=\<f,\phi_j\> $, 
where   $\phi_1,\dots,\phi_m\in \cH$.  Therefore,  $M(f)=(\<f,\phi_1\>,\dots,\<f,\phi_m\>)\in \R^m$.  
We let $W:=\span\{\phi_j\}_{j=1}^m$, which is an $m$ dimensional subspace of $\cH$.  
 We can always perform a Gram-Schmidt orthogonalization and assume therefore    that 
 $\phi_1,\dots,\phi_m\in \cH$ is an orthonormal basis for $W$ (see Remark \ref{equivinfo}).   
 We have  $\cN=W^\perp$. From (\ref{Quotdual2}) and (\ref{Quotdual1}) we infer that $\|\cdot \|_M $ 
 on $\R^m$  is the $\ell_2(\R^m)$ norm.   Therefore, the approximation  map is simple 
 least squares fitting.  Namely, to our data $w$, we find the element $ z^*(w)\in Z$, where $Z:=M(V)$, such that
$$z^*(w):=\argmin_{z\in Z} \sum_{j=1}^m |w_j-z_j|^2.
$$
The element $v^*(w)=M_V^{-1}(z^*(w))$ is the standard least squares fit to the data $(f(P_1),\dots, f(P_m))$ by 
vectors $(v(P_1),\dots,v(P_m))$ with $v\in V$, and is found by the usual matrix inversion in least squares.   
This gives the best approximation to $w$ in $\|\cdot\|_M$ by the elements of $Z$,  and hence $\lambda=1$.
 The lifting  $\Delta(w_1,\ldots,w_m):=\sum_{j=1}^m w_j \phi_j$ is linear and $\|\Delta\|=1$.  Hence, we have the algorithm
 
 \be
 \label{l2algorithm}
 A(w) = M_V^{-1}(z^*(w))+\Delta(w-z^*(w))= v^*(w)+\sum_{j=1}^m [w_j-z_j^* (w)]\phi_j,
 \ee
 which is the algorithm presented in \cite{MPPY} and further studied in \cite{BCDDPW1}.
 The sum in \eref{l2algorithm} is a correction so that $A(w)\in\cK_w$, i.e., $M(A(w))=w$.
   
   \begin{remark}   Our general theory says that the above algorithm is near  optimal  
   with constant $2$ for recovering $\cK_w$.  It is shown in {\rm \cite{BCDDPW}} that, in this case,  
   it is actually an  optimal algorithm.   The reason for this is that   the sets $\cK_w$ in this
    Hilbert case setting have a center of symmetry, so  Proposition {\rm \ref{prop:centrallysymmetric}}  can be applied.  \end{remark}
   
\begin{remark}  It was shown in {\rm \cite{BCDDPW1}}  that the calculation can be streamlined   
by choosing at the beginning certain favorable  bases  for $V$ and $W$.   In particular, the quantity $\mu(\cN,V)$
can be  immediately computed from the cross-Grammian of the favorable bases.
\end{remark} 

\subsection{Example 1}
\label{ss:ex1}  In this section, we  summarize how the algorithm works for Example 1.  Given 
$P_j\in D$, $j=1,\ldots,m$, $P_i\neq P_j$,
and the data $w=M(f)=(f(P_1),\cdots,f(P_m))$,  the first step is to find
the min-max approximation to $w$ from the space $Z:=M(V)\subset \R^m$. In other words, we find 
\be
\label{minmax}
 z^*(w):=\argmin_{z\in Z} \max_{1\le j\le m}|f(P_i)-z_i|=\argmin_{v\in V} \max_{1\le j\le m}|f(P_i)-v(P_i)|.
\ee
Note that for general  $M(V)$ the point $z^*(w)$ is not necessarily unique.  For certain $V$, however, we have  uniqueness.  

Let us consider the case when $D=[0,1]$ and $V$ is a Chebyshev space on $D$, i.e., for any $n$ points $Q_1,\dots,Q_n\in D$, and any data $y_1,\dots,y_n$, there is a unique function $v\in V$ which satisfies $v(Q_i)=y_i$, $1\le i\le n$. In this case, when $m=n$, problem \eref{minmax} has a unique solution 
$$
z^*(w)=w=M(v^*(w))=(v^*(P_1), \ldots,v^*(P_m)),
$$
where $v^*\in V$ is the unique interpolant to the data $(f(P_1),\cdots,f(P_m))$ at the points $P_1, \ldots,P_m$.
 For $m\geq n+1$, let us denote by $V_m$ the restriction of $V$ to the point set  $\Omega:=\{P_1,\ldots,P_m\}$. Clearly, $V_m$ 
 is a Chebyshev space  on $C(\Omega)$ as well,   and therefore 
there is a unique   point  $z^*(w):=(\tilde v(P_1), \ldots,\tilde v(P_m))\in V_m$,   coming from the evaluation of a unique 
$\tilde v\in V$, which is the  best approximant from $V_m$ to $f$ on $\Omega$.
The point $z^*(w)$ is characterized by an oscillation property.   Various algorithms for finding $\tilde v$ are known and go under the name
Remez algorithms.

 In  the general case where $V$ is not necessarily a Chebyshev space,  a minimizer $z^*(w)$  can still be  found by convex minimization, and the approximation 
mapping $\Lambda$ maps $w$ to a $z^*(w)$.  Moreover,  $z^*(w)=M(v^*(w))$ for 
some $v^*(w)\in V$, where $v^*(w) $ is characterized by solving the minimization
$$
v^*(w)=\argmin_{v\in V}\|w-M(v)\|_M= \argmin_{v\in V}\inf_{g:\,M(g)=w}\|g-v\|_{\cX}=\argmin_{v\in V} \dist(v,\cX_w).
$$

We have seen that the lifting in this case is simple.   We may take functions $\psi_j\in C(D)$, with disjoint supports and of norm one, such that $\psi_i(P_j)=\delta_{i,j}$.   Then, we can take our lifting to be the operator that maps $w\in\R^m$ into the function $\sum_{j=1}^mw_j\psi_j$.   This is a linear lifting with norm one.  Then, the  algorithm $A$ is given by 
\be
\label{optalgc}
A(w):= M_V^{-1}(z^*(w))+\sum_{j=1}^m(w_j-z_j^*(w))\psi_j=v^*(w)+\sum_{j=1}^m(w_j-z_j^*(w))\psi_j,\quad w\in\R^m.
\ee
The sum in \eref{optalgc} is a correction to $v^*(w)$ to satisfy the data. From \eref{tnb}, we know  that 
for each $w\in\R^m$, we have
$$
 \sup_{f\in\cK_w} \|f-A(w)\|\le   2\rad(\cK_w),
$$
and so the algorithm is near  optimal with constant $2$ for each of the classes $\cK_w$.   

To give an a priori bound for the  performance of this algorithm, we need to compute     $\mu(\cN,V)$.

 \begin{lemma}
\label{peval}
Let ${\cal X}=C(D)$, $V$ be a subspace of $C(D)$,   and $M(f)=(f(P_1),\ldots,f(P_m))$, 
where $P_j\in D$, $j=1,\ldots,m$ are $m$ distinct points in $D\subset\R^d$.
Then, for  $\cN$  the null space of $M$, we have
$$
 \frac{1}{2}\sup_{v\in V}\frac{\|v\|_{C(D)}}
 {\displaystyle{\max_{1\le j\le m}|v(P_j)}|}\leq\mu({\cal N}, V)\leq 2 \sup_{v\in V}\frac{\|v\|_{C(D)}}{\displaystyle{\max_{1\le j\le m}|v(P_j)|}}.
$$

\end{lemma}
\noindent
{\bf Proof:} From Lemma \ref{mu} and Lemma \ref{comp}, we have
\be
\label{tu}
\frac{1}{2}\|M_V^{-1}\|\le \mu({\cal N}, V)\leq 2\|M_V^{-1}\|.
\ee
Since, we know $\|w\|_M=\max_{1\le j\le m}|w_j|$, we obtain that 
$$\|M_V^{-1}\|= \sup_{v\in V}\frac{\|v\|_{C(D)}}{\displaystyle{\max_{1\le j\le m}|v(P_j)}|},
$$ 
and the lemma follows. 
\hfill $\blacksquare$

  From \eref{algper}, we obtain the a priori performance bound %
 \be
  \label{globalnearbest1}
\sup_{w\in\R^m} \sup_{f\in\cK_w} \|f-A(w)\|_{C(D)}\le   4\e\mu(\cN,V).
\ee
Moreover, we know from Theorem \ref{theorem:main} that \eref{globalnearbest1} cannot be improved by any algorithm except for the possible 
removal of the factor $4$,  and hence the algorithm is globally near optimal.

\begin{remark}  
\label{rem:io}  
 It is important to note that the algorithm $A:\ w\rightarrow A(w)$ does not depend on $\e$, and so one obtains 
 for any $f$ with the data $w=(f(P_1),\dots,f(P_m))$ the performance bound
$$ 
 \|f-A(w)\|_{C(D)}\le 4\mu(\cN,V)\dist(f,V).
$$
Approximations of this form are said to be instance optimal with constant $4\mu(\cN,V)$.  
\end{remark}

As an illustrative example,  consider the space $V$ of trigonometric polynomials of degree $\le n$ on 
$D:=[-\pi,\pi]$, which is a Chebyshev system  of dimension $2n+1$. 
We take $\cX$ to be the space of continuous functions on  $D$ which are periodic, 
i.e.,  $ f(-\pi)=f(\pi)$.  If the data consists of the values of $f$ at $2n+1$  distinct points $\{P_i\}$,
then the min-max approximation is simply 
the interpolation projection  $\cP_nf$ of  $f$  at these points and $A(M(f)))=\cP_nf$. 
The  error estimate for this case is
$$
\|f- \cP_nf\|_{C([-\pi,\pi])}\le (1+\|\cP_n\|)\dist(f,V).
$$
It is well known (see \cite{Z}, Chapter 1 of Vol. 2)  that for $P_j:= -\pi+j\frac{2\pi}{2n+1}$, $j=1,\dots, 2n+1$, 
$\|\cP_n\|\approx \log n$.  However,
if we double the number of points, and keep them equally spaced,  then     it is  known that $\|M_V^{-1}\|\le 2$ (see  \cite{Z}, Theorem 7.28).  Therefore from  \eref{tu},  we obtain 
 $\mu(\cN,V)\leq 4$, and we derive the bound
\be
\label{btrig}
\|f-A(M(f)))\|_{C([-\pi,\pi])}\le  16\dist(f,V).
\ee

\subsection{Example 2} 
\label{ss:ex2} This case is quite similar to Example 1.  The main difference is that now
\be
\label{lpmin}
z^*(w):=\argmin_{z\in Z} \|w-z\|_{\ell_p(\R^m)},
\ee
and hence when $1<p<\infty$ it can be  found by minimization in a uniformly convex norm.  
We can take the lifting $\Delta$ to be $\Delta(w)=\sum_{j=1}^m w_j\psi_j$, where now 
$\psi_j$ has the same support as $g_j$  and  
$L_{p}(D)$ norm one, $j=1,\ldots,m$.   The algorithm is again given by
\eref{optalgc}, and is near  optimal  with constant $2$ on each class $\cK_w$, $w\in\R^m$, that is
$$
\|f-A(M(f))\|_{L_p(D)}\le 2\rad(\cK_w) \le 4\mu(\cN,V)\e ,
$$
where the last inequality follows from \eref{algper}.
 
   Similar to Lemma \ref{peval}, we have  the following bounds for $\mu(\cN,V)$,
   $$
 \frac{1}{2}\|M_V^{-1}\|\le \mu({\cal N}, V)\leq 2\|M_V^{-1}\|,
$$
where now the norm of $M_V^{-1}$ is taken as the operator norm from $L_p(D)$ to $\ell_p(\R^m)$, and hence is 
$$
\|M_V^{-1}\|= \sup_{v\in V}\frac{\|v\|_{L_p(D)}}{\|(l_1(v),\ldots,l_m(v))\|_{\ell_p(\R^m)}}.
$$ 
 
  
\subsection{Example 3} 
\label{ss:ex3} As mentioned earlier, our interest in Example 3 is because it illustrates certain theoretical features.   
In this example, the norm $\|\cdot\|_M$ is the $\ell_\infty(\R^m)$ norm,  and approximation in this norm
was already discussed in Example 1.   The interesting aspect of this example centers around liftings.  
We know that any linear lifting must have norm $\geq\sqrt{m/2}$.   On the other hand, 
we have given in \eref{Riesz1} an explicit formula for a (nonlinear) lifting with norm one.  
So, using this lifting, the algorithm $A$ given in \eref{algmap}  will be near  optimal 
with constant $2$ for each of the classes $\cK_w$.  

\section{Relation to sampling theory}  
\label{sec:sampling}
 The results we have put forward, when restricted to problems of sampling,  have some overlap with recent results.  In this section, we point out
these connections and what new light our general theory sheds on sampling problems.  The main point to be made is that our results give a general framework for sampling in Banach spaces that includes  many of the specific examples studied in the literature.

 Suppose that $\cX$ is a Banach space and $l_1,l_2,\dots, $ is a possibly infinite sequence of
linear functionals from $\cX^*$.  The application of the $l_j$ to an $f\in \cX$ give a sequence of samples of $f$.  Two prominent examples are the following. 
\vskip .1in
\noindent
{\bf  Example PS: Point samples of continuous functions.}  Consider the space $\cX=C(D)$ for a domain $D\subset\R^d$ and a sequence of points $P_j$ from $D$.  Then, the point evaluation functionals $l_j(f)=f(P_j)$, $j=1,2,\dots$, are point samples of $f$.   Given a compact subset $K\subset \cX$, we are interested in how well we can recover $f\in K$ from the information $l_j(f)$, $j=1,2,\dots$.

\vskip .1in 
\noindent
{\bf Example FS: Fourier samples.}  Consider the space $\cX=L_2(\Omega)$, $\Omega=[-\pi,\pi]$, and the linear functionals   
\be
\label{fouriercoefficients}
l_j(f):= \frac{1}{2\pi}\int_\Omega f(t)e^{-ijt}\, dt,\quad j\in \Z,
\ee
which give the Fourier coefficients of $f$.
 Given a compact subset $K\subset \cX$, we are interested in how well we can recover $f\in K$ in the norm of $L_2(\Omega)$ 
 from the information $l_j(f)$, $j=1,2,\dots$.
 \vskip .1in
 The main problem in sampling is to build  reconstruction operators $A_m:\R^m\mapsto \cX$ such that  the reconstruction mapping $R_m(x):=A_m(M_m(x))$ provide
 a good approximation to $x$. Typical questions are (i) Do there exist such mappings such that $R_m(x)$ converges to $x$ as $m\to \infty$,  for each $x\in\cX$?, (ii) What is the best performance in terms of rate of approximation on specific  compact sets $K$?, (iii) Can we guarantee the stability of these maps in the sense of how they perform with respect to noisy observations?.
 
 The key in connecting such sampling problems with our theory is that the compact sets $K$ typically considered are either directly defined by   approximation or can be equivalently
 described by such approximation.  That is, associated to $K$ is a sequence of spaces $V_n$, $n\ge 0$, each of dimension $n$, and $f\in K$ is equivalent to
 \be
 \label{as1}
 \dist(f,V_n)_X\le \e_n,\quad n\ge 0,
 \ee
  where $(\e_n)$ is a known sequence of positive numbers which decrease to zero.
 Typically, the $V_n$ are nested, i.e. $V_n\subset V_{n+1}$, $n\ge 0$.
 Such characterizations of sets $K$ are often provided by the theory of approximation.  For example, a periodic function $f\in C[-\pi,\pi]$ is in Lip $\alpha$, $0<\alpha<1$ if and only if
 \be
 \label{lipspaces}
 \dist(f,\cT_n)_{C[-\pi,\pi]}\le C(f)(n+1)^{-\alpha}, \quad n\ge 0,
 \ee
 with $\cT_n$ the space of trigonometric polynomials of degree at most  $n$ and moreover, the Lip $\alpha$ semi-norm is equivalent to the smallest constant $C(f)$ for which
 \eref{lipspaces} holds.  Similarly, a function $f$ defined on $[-1,1]$ has an analytic extension to the region in the plane with boundary given by Bernstein ellipse $E_\rho$ if and only if 
 \be
 \label{Eanalytic}
 \dist(f,\cP_n)_{C[-1,1]}\le C(f)\rho^{-n},\quad n\ge 0,
 \ee
 where $\cP_n$ is the space of algebraic polynomials of degree at most  $n$ in one variable (see \cite{PTK}).
 
For the remainder of this section, we assume that the set $K$ is 
\be
\label{defK}
K:=\{ x\in\cX: \dist(x,V_n)_\cX\le \e_n\}=\bigcap_{n\ge 0}\cK(\e_n,V_n).
\ee
       Our results previous to this section assumed only the knowledge that $f\in \cK(\e_n,V_n)$ for one fixed value of $n\le m$.  

Our general theory (see Theorem \ref{theorem:ba}) says that given the first $m$ samples 
\be
\label{sample}
M_m(x):=(l_1(x),\dots,l_m(x)),\quad x\in\cX,
\ee
 then for any $n\le m$, the mapping $A_{n,m}$ from $\R^m\mapsto \cX$, given by \eref{algmap},   provides an approximation $A_{n,m}(M_m(x))$ to $x\in K$ with the accuracy
 \be
 \label{sample1}
 \|x-A_{n,m}(M_m(x))\|_\cX\le C\mu(V_n,\cN_m)\e_n,
 \ee
 where $\cN_m$ is the null space of the mapping $M_m$.  Here, we know that theoretically $C$ can be chosen as close to $8$ as we wish but in numerical implementations, depending on the specific setting, $C$ will generally be a known constant but larger than $8$.  In the two above examples,  one can take $C=8$ both theoretically and numerically.
 \begin{remark}
 \label{switch}
 It is more convenient in this section to use the quantity $\mu(V,\cN)$ rather than $\mu(\cN,V)$.  Recall that $\frac{1}{2}\mu(V,\cN)\le \mu( \cN,V))\le 2\mu(V,\cN)$ and therefore
 this switch only effects constants by a factor of at most $2$.
 \end{remark}
 \begin{remark}
 \label{samplingremark}
  Let  $A_{n,m}^*:\R^m\mapsto V_n$ be the mapping defined by \eref{algmap} with the second term   $\Delta(w-M(w))$ on the right deleted.  Fom \eref{equiv}, it follows that the term that is dropped has norm not exceeding $\|\Delta\|\e_n$ whenever $x\in K$, and since we can
  take $\|\Delta\|$ as close to one as we wish, the resulting operators satisfy \eref{sample1}, with a new value of $C$,  but now they map into $V_n$.
  \end{remark}

  Given a value of $m$ and the information map $M_m$, we  are allowed to choose $n$, i.e., we can choose the space $V_n$.   Since $\e_n$ is known, from the point of view of the error bound \eref{sample1}, given the $m$ samples, the best choice of $n$ is   
 \be
 \label{sample2}
 n(m):=\argmin_{0\le n\le m} \mu(V_n,\cN_m)\e_n,
 \ee
 and gives the bound
  \be
 \label{sample11}
 \|x-A_{n(m),m}(M_m(x))\|_\cX\le C \min_{0\le n\le m}\mu(V_n,\cN_m)\e_n,\quad x\in K.
 \ee

 This brings out the importance of giving good estimates for $\mu(V_n,\cN_m)$ in order to be able to select the best choice for $n(m)$.   Consider the case of point samples.   Then, the results of Example 2 in the previous section show that in this case, we have
 \be
 \label{sample3}
   \mu(V_n,\cN_m)= \sup_{v\in V_n}\frac{\|v\|_{C(D)}}{\displaystyle{\max_{1\le j\le m}|v(P_j)}|}.
 \ee
 The importance of this ratio of the continuous and discrete norms, in the case $V_n$ are spaces of algebraic polynomials,  has been known for some time.  It is equivalent to the Lebesgue constant when $n=m$
 and has been recognized as an important ingredient in sampling theory dating back at least to Sch\"onhage \cite{Sc}.   
 
 A similar ratio of continuous to discrete norms  determine $\mu(V_n,\cN_m)$ for other sampling settings.  For example, in the case of   the  Fourier sampling ({\bf Example FS} above), we have for any space $V_n$ of dimension $n$
 \be
 \label{FS1}
 \mu(V_n,\cN_{2m+1})= \sup_{v\in V_n}\{\|v\|_{L_2(\Omega)}:\|\cF_m(v)\|_{L_2(\Omega)}=1\},
 \ee
 where $\cF_m(v)$ is the $m$-th partial sum of the Fourier series of $v$.   The right side of \eref{FS1}, in the case $V_n=\cP_{n-1}$  is studied in \cite{AHS} (where it is denoted by  $B_{n,m}$).
  Giving bounds for  quotients, analogous to  those   in \eref{sample3}, has been  a central topic in sampling theory (see \cite{AH,AHP, AHS, PTK}) and such bounds have been obtained in specific settings, such as the case of equally spaced point samples on $[-1,1]$ or Fourier samples.  The present paper does not contribute to the problem of estimating
  such ratios of continuous to discrete norms.

 The results of the present paper give a general framework for the  analysis of sampling.   Our construction of the operators $A_{n,m}$ (or their modification $A_{n,m}^*$ given by Remark \ref{samplingremark}),  give performance bounds that match those given in the literature in specific settings such as the two examples given at the beginning of this section.
 It is interesting to ask in what sense these bounds are optimal.  Theorem \ref{theorem:ba} proves optimality of the bound \eref{sample1} if the assumption that $f\in K$ is replaced by the less demanding assumption that $f\in \cK(\e_n,V_n)$, for this one fixed value of $n$.
 The knowledge
 that   $K$ satisfies $\dist(K,V_n)\le \e_n$, for all $n\ge 0$, could allow an improved performance, since it is a more demanding assumption.  In the case of a Hilbert space, this was shown   to be the case in \cite{BCDDPW} where, in some instances,  much improved bounds were obtained from this additional knowledge.  However, the examples in \cite{BCDDPW}  are not for classical settings such as polynomial or trigonometric polynomial approximation.   In these cases, there is no known improvement over the estimate \eref{sample11}.  
 
Next, let us consider  the question of whether, in the case of an infinite number of samples, the samples guarantee that every $x\in K$ can be approximated to arbitrary accuracy from these samples.   This is of course the case whenever
 \be
 \label{sample4} \mu(V_{n(m)},\cN_m)\e_{n(m)}\to 0, \quad m\to\infty.
 \ee
 In particular, this will be the case whenever the sampling satisfies that for each finite dimensional space
 \be
 \label{fds}
 \lim_{m\to\infty}\mu(V,\cN_m)\le C,
 \ee 
 for a fixed constant $C>0$, independent of $V$.  Notice that the spaces $\cN_m$ satisfy $\cN_{m+1}\subset \cN_m$ and hence the sequence $(\mu(V,\cN_m))$ is non-increasing.
 
 Of course a necessary condition for \eref{fds} to hold is that the sequence of functionals $l_j$, $j\ge 1$, is total on $\cX$, i.e., for each $x\in\cX$, we have
 \be
 \label{total}
 l_j(x)=0, \ j\ge 1,\implies x=0.
 \ee
 It is easy to check that when the $(l_j)_{j\ge 1}$ are total, then for each $V$, we have that \eref{fds} holds with a constant $C_V$ depending on $V$.  
  Indeed, if \eref{fds} fails  then $\mu(V,\cN_m)=+\infty$ for all $m$ so by our previous remark $V\cap\cN_m\neq\{0\}$.  The  linear spaces $ V\cap \cN_m$, $m\ge 1$, are nested and contained in $V$.     If $V\cap \cN_m\neq \{0\}$  for all $m
  $,    then there is a $v\neq 0$ and $v\in \bigcap_{m\geq 1} \cN_m$.  This  contradicts \eref{total}. 
  
   However, our interest is to have a bound uniform in $V$.     To derive such a uniform bound, we introduce the following notation.  Given the sequence $(l_j)_{j\ge 1}$, we let  $L_m:=\span\ \{l_j\}_{j\leq m} $ and $L:=\ \span\{l_j\}_{j\ge 1}$ which are  closed  linear subspaces  of $\cX^*$.  We denote by $U(L_m)$ and $U(L)$ the unit ball of these spaces with the $\cX^*$ norm.
 For any $0<\gamma \leq 1$, we say that the  sequence $(l_j)_{j\ge 1}$ is $\gamma$-norming,     if   we have
 \be
 \label{kappanorming}
  \sup_{l\in U(L)}|l(x)|\geq \gamma \|x\|_\cX,\quad x\in\cX.
  \ee
   Clearly any $\gamma$-norming sequence is total. If $\cX$ is a reflexive Banach space, then, the Hahn-Banach theorem gives that  every total sequence $(l_j)$  is $1$-norming.

 \begin{theorem}
 \label{theorem:total}  If  $\cX$ be any  Banach space and  suppose that the functionals $l_j$, $j=1,2,\dots$, are  $\kappa$-norming, then for any finite dimensional subspace $V\subset \cX$,
 we have  
 \be
 \label{tt1}
\lim_{m\to \infty} \mu(V,\cN_m) \leq \gamma^{-1}.
 \ee
 \end{theorem}
 
 \noindent 
 {\bf Proof:}  Since the unit ball $U(V)$ of $V$ in $\cX$ is compact, for any  $\delta>0$,  there exists an $m$ such that
 \be
 \label{uniV}
 \sup_{l\in U(L_m)} |l(v)|\geq \frac{\gamma}{1+\delta}\|v\|_\cX,\quad v\in V.
 \ee
 If we fix this value of $m$, then  for any $v\in U(V)$   and any  $\eta\in \cN_m$ we have
 $$\|v-\eta\|_\cX\geq \sup_{l\in U(L_m)}
  |l(v-\eta)|= \sup_{l\in U(L_m)}
  |l(v)|   \ge  \frac{\gamma}{1+\delta}. $$
 From \eref{defmu} follows that $\mu(V,\cN_m)\leq (1+\delta)\gamma^{-1}$ and since $\delta $ was arbitrary this proves the theorem.\hfill $\Box$

 \begin{cor}
\label{correcover}
Let  $\cX$ be any separable   Banach space and let $(l_j)_{j\ge 1}$  be any sequence of  functionals from $\cX^*$ which are $\gamma$ norming for some $0<\gamma\le 1$.   Then, we have the following results:

\vskip .1in
\noindent
{\rm (i)} If  $(V_n)_{n\ge 0}$  is any sequence of nested finite dimensional spaces whose closure   is  $\cX$, then, for each $m$, there is  a choice $n(m)$ such that
 \be
\label{satisfy1}
\|x-A^*_{n(m),m}(M_m(x))\|_\cX\le C\mu(V_{n(m)},\cN_m)\dist(x,V_{n(m)})_\cX \le 2C\dist(x,V_{n(m)})_\cX,
\ee 
with $C$ an absolute constant, and the right side of \eref{satisfy1} tends to zero as $m\to\infty$. 

\vskip .in
\noindent
{\rm (ii)} There exist operators $A_m$ mapping $\R^m$ to $\cX$ such that $A_m(M_m(x))$ converges to $x$, for all $x\in \cX$.
\vskip .1in
\noindent
In particular, both {\rm (i)} and {\rm (ii)} hold whenever $\cX$ is reflexive and $(l_j)_{j\ge 1}$ is  total.
\end{cor}

\noindent
{\bf Proof:}  In view of the previous theorem,  for each sufficiently large $m$, there is a $n(m)\ge 0$ such that $\mu(V_{n(m)},\cN_m)\le 2\gamma^{-1}$ and we can take $n(m)\to \infty$ as $m\to \infty$.  Then, (i) follows  from \eref{sample11}.  The statement (ii) follows from (i) because there is always a sequence $(V_n)$ of spaces of dimension $n$ whose closure is dense in $\cX$.
\hfill $\Box$

While the spaces $C$ and  $L_1$,    are not reflexive, our two examples are still covered by Theorem \ref{theorem:total} and Corollary \ref{correcover} 
\vskip .1in

\noindent
{\bf Recovery for Example PS:}  If the points $P_j$ are dense in $C(D)$, then the sequence of functionals  $l_j(f)=f(P_j)$, $j\ge 1$ is $1$ norming and Theorem \ref{theorem:total} and Corollary \ref{correcover} hold for this sampling.

\vskip .1in
\noindent
{\bf Recovery for Example FS:}  The sequence $(l_j)_{j\ge 1}$ of Fourier samples is $1$ norming  for each of the spaces $L_p(\Omega)$, $1\le p<\infty$, or $C(\Omega)$  and hence Corollary \ref{correcover} hold for this sampling. 

\vskip .1in
\noindent
We leave the simple details of these last two statements to the reader.

 Let us note that, in general,  the totality of the linear functionals is not sufficient to guarantee Theorem \ref{theorem:total}.  A simple example
 is to take $\cX=\ell_1=\ell_1(\N)$, with its usual basis $(e_j)_{j=1}^\infty$ and coordinate functionals   $(e_j^*)_{j=1}^\infty$. We fix any  number $a>1$ and consider the linear  functionals
 \begin{eqnarray*}
 l_1&:=& ae_1^*-\sum_{j=2}^\infty e_j^*, \\
 l_k&:=&e_k^* , \quad  k\ge 2.
 \end{eqnarray*}
We then have
 \begin{enumerate}
 \item {\em The system $(l_j)_{j=1}^\infty$ is total.} That is,  if for $x\in \cX$ we have  $l_j(x)=0$ for $j=1,2,\dots$ then $x=0$.  Indeed, if  $x=(x_j)_{j\ge 1}\in \ell_1$ with $l_k(x)=0$ for $k\geq 2$, then $x=\mu e_1$.  The requirement  $l_1(x)=0$ then implies that $\mu=0$ and hence $x=0$.
 \item {\em For $m\ge 1$, $\cN_m$ equals the set of all vectors  $(\eta_1,\eta_2,\dots)\in \ell_1$ such that $\eta_1=a^{-1}\sum_{j> m} \eta_j$ and $\eta_2=\cdots = \eta_{m}=0$.} 
  \item {\em $\dist (e_1,\cN_m)_{\ell_1}=a^{-1}$, $m\geq 1$.} Indeed, if  $\eta\in \cN_m$,  then by (ii),  $\eta=(\eta_1,0,\cdots, 0,\eta_{m+1},\cdots)$ and $\eta_1=a^{-1}\sum_{j>m}\eta_ j=:a^{-1}z$.   Therefore,    we have 
 $$
 \|e_1-\eta\|_{\ell_1}= |1-a^{-1}z|+\sum_{j>m}|\eta_j| \ge |1-a^{-1}z|+|z|\ge |1-a^{-1}|z||+|z|
 $$ 
 and therefore, we have  $\displaystyle{\inf_{\eta\in \cN_m}\|e_1-\eta\|=\inf_{z>0}|1-az|+z=a^{-1}}$.
 \end{enumerate}
  It follows from these remarks that for  $V$  the one dimensional space spanned by  $e_1$, we have $\mu(V,\cN_m)\ge  a$, for all $m\ge 1$.  Let us note that in fact this system of functionals is $a^{-1}$ norming.

One can build on this example to construct  a non-reflexive space $\cX$ and a sequence of functionals $l_j\in \cX^*$ which are total but for each $j\ge 1$ there is a $V_j$ of dimension one
such that
\be
\label{Vj}
\lim_{m\to\infty} \mu(V_j,\cN_m)\ge j.
\ee

 A major issue in generalized sampling is whether the reconstruction maps are numerically stable.  To quantify stability, one introduces for any potential reconstruction maps
 $\phi_m:\R^m\mapsto \cX$, the operators $R_m(f):=\phi_m(M_m(f))$ and the condition numbers
 \be
 \label{defcond}
 \kappa_m:= \sup_{f\in\cX}\lim_{\e\to 0} \sup_{\|g\|_\cX=1}\frac{\|R_m(f+\e g)-R_m(g)\|_\cX}{\e},
 \ee
 and asks whether there is a uniform bound on the $\kappa_m$.  In our case  
 \be
 \label{phim}
  \phi_m=A_{n(m),m}^* = M_{V_{n(m)}}^{-1}\circ \Lambda_m,  
  \ee
 where $\Lambda_m$ is the approximation operator for approximating $M_m(f)$ by the elements from $M_m(V_{n(m)})$.  Thus,
 \be
 \label{cond1}
 \kappa_m\le  \|M_{V_{n(m)}}^{-1}\|_{{\rm Lip}\  1}\|\Lambda_m\|_{{\rm Lip} \ 1 }\le C\mu(V_{n(m)},\cN_m)\|\Lambda_m\|_{{\rm Lip} \ 1 },
 \ee
 where $\|\cdot\|_{{\rm Lip}\ 1}$ is the Lipschitz norm of the operator.   Here to bound the Lipschitz norm of the operator $M_{V_{n(m)}}^{-1}$ we used the fact that it
 is a linear operator and hence its Lipschitz norm is the same as its norm and in this case is given by  (i) of Lemma \ref{comp}.  Under the assumptions of Theorem \ref{theorem:total}, we know
 that $\mu(V_{n(m)},\cN_m)\le C$ for a fixed constant $C$.

In the case $\cX$ is a Hilbert space, the approximation operator $\Lambda_m$ can also be taken as a linear operator of norm one on a Hilbert space and hence we obtain
a uniform bound on the condition numbers $\kappa_m$.  For more general Banach spaces $\cX$, bounding the Lipschitz norm of $\Lambda_m$ depends very much on the specific 
spaces $V_n$, $\cX$, and the choice of $\Lambda_m$.   However, from the Kadec-Snobar theorem, we can always take for $\Lambda_m$ a linear projector whose norm is bounded by $\sqrt{n(m)}$ and therefore,  under the assumptions of Theorem \ref{theorem:total}, we have
\be
\label{boundkappa}
\kappa_m\le C \sqrt{n(m)},\quad m\ge 1.
\ee

\section{Choosing measurements} 
\label{sec:choosemeasurements}  
In some settings, one knows the space $V$,  but is allowed to choose  the measurement functionals $l_j$, $j=1,\dots,m$.  
In this section, we discuss how our results can be a guide in such a selection.   The main issue is to keep $\mu(\cN,V)$ as small as possible
for this choice,  and so we concentrate on this.

\noindent
Let us recall that from Lemma \ref{mu} and Lemma \ref{comp}, we have  
$$
 \frac{1}{2}\|M_V^{-1}\|\le \mu({\cal N}, V)\leq 2\|M_V^{-1}\|.
$$
Therefore, we want to choose $M$ so as to keep 
$$
\|M_V^{-1}\|= \sup_{v\in V}\frac{\|v\|_\cX}{\|M_V(v)\|_M}
$$
small.   In other words, we want to keep $\|M_V(v)\|_M $ large whenever $\|v\|_\cX=1$.

\noindent
{\bf Case 1:} Let us first consider the case when $m=n$.   Given any linear functionals $l_1,\dots,l_n$,  
which are linearly independent over $V$ (our candidates for measurements), we can choose a basis for $V$ 
which is dual to the $l_j$'s,  that is, we can choose $\psi_j\in V$, $j=1,\dots,n$, such that
$$
l_i(\psi_j)=\delta_{i,j},\quad 1\le i,j\le n.
$$
It follows that each $v\in V$ can be represented as $v=\sum_{j=1}^nl_j(v)\psi_j$.  
The operator  $P_V: \cX\rightarrow V$, defined as 
\be
\label{projV}
P_V(f)= \sum_{j=1}^n l_j(f)\psi_j,\quad f\in\cX,
\ee
is a projector from $\cX$ onto $V$, and any projector onto $V$ is of this form.  If we take  
$M(v)=(l_1(v),\dots,l_n(v))$, we have
\be
\label{keep1}
\|M(v)\|_M=  \inf_{M(f)=M(v)}\|f\|=  \inf_{P_V(f)=v}\|f\|=\inf_{P_V(f)=v} \frac{\|f\|}{\|P_V(f)\|} \|v\|.
\ee
If we now take the infimum over all $v\in V$ in \eref{keep1}, we run through all $f\in \cX$, and hence
$$
\inf_{\|v\|=1}\|M(v)\|_M= \inf_{f\in\cX}  \frac{\|f\|}{\|P_V(f)\|} =\|P_V\|^{-1}.
$$
In other words,
$$
\|M_V^{-1}\| =\|P_V\|.
$$
This means the best choice of measurement functionals is to take the linear projection onto $V$ with smallest norm, 
then take any basis $\psi_1,\dots,\psi_n$ for $V$ and represent the projection in terms of this basis as in \eref{projV}.   
The dual functionals  in this representation are  the measurement functionals.

Finding projections of minimal norm onto a given subspace $V$ of a Banach space $\cX$ is a well-studied problem in functional analysis.
A famous theorem of Kadec-Snobar \cite{KS} says that there always exists such a projection with
\be
\label{ks}
\|P_V\|\le \sqrt{n}.
\ee
It is known that there exists Banach spaces $\cX$ and subspaces $V$ of dimension $n$, where \eref{ks}  cannot be improved in 
the sense that for any projection onto $V$ we have $\|P_V\|\ge c\sqrt{n}$ with an absolute constant $c>0$.  If we translate this result to our setting
of recovery, we see that given $V$ and $\cX$ we can always choose measurement functionals $l_1,\dots,l_n$,  such that $\mu(\cN,V)\le 2\sqrt{n}$, 
and this is the best we can say in general.   

\begin{remark}
\label{rem:ks}
For a general Banach space $\cX$ and a finite dimensional subspace $V\subset \cX$ of dimension $n$, finding a minimal norm projection 
or even a near minimal norm projection onto $V$ is not constructive.   There are related procedures such as Auerbach's 
theorem {\rm\cite[II.E.11]{Wbook}},   which give the poorer estimate $Cn$ for the norm of $\|P_V\|$.   These constructions are easier to describe 
but they also are not computationally feasible.  
\end{remark}

 \begin{remark}
 \label{rem:ks1}
 If $\cX$ is an $L_p$ space, $1<p<\infty$, then the best bound in \eref{ks} can be replaced by $n^{|1/2-1/p|}$, and this is again known to be optimal, 
 save for multiplicative constants. When $p=1$ or $p=\infty$ (corresponding to $C(D)$), we obtain the best bound $\sqrt{n}$ and this cannot 
 be improved for general $V$. Of course, for specific $V$ the situation may be much better.  Consider $\cX=L_p([-1,1])$,  and $V=\cP_{n-1}$ 
 the space of polynomials of degree at most $n-1$.  In this case, there are projections with norm $C_p$,  depending only on $p$. For example,
 the projection given by the Legendre polynomial expansion has this property.  For $\cX=C([-1,1])$, the projection given by interpolation
 at the zeros of the Chebyshev polynomial of first kind has norm $C\log n$,  and this is again optimal save for the constant $C$.
\end{remark}

\noindent
{\bf Case 2:} 
Now, consider the case when the number  of measurement functionals $m>n$.  
One may think that one can drastically improve
on the results for $m=n$.  We have already remarked that this is possible in some settings 
by simply doubling the number of data functionals (see \eref{btrig}).  
While adding additional measurement functionals does decrease $\mu$,  generally speaking, 
we must  have $m$ exponential in $n$ to guarantee that $\mu$ is independent of $n$. 
To see this,  let us discuss  one special case of Example 1. We fix 
$D=:\{x\in \R^n\, :\ \sum_{j=1}^n x_j^2 = 1\}$ and the  subspace $V\subset C(D)$ of all linear 
functions restricted to $D$,  i.e.,  $f\in V$ if and only if 
$$
f(x)=f_a(x):=\sum_{j=1}^n a_jx_j, \quad a:=(a_1,\ldots,a_n)\in \R^n.
$$
 It is obvious that $V$ is an $n$ dimensional subspace. Since for $f\in V$,  we have $\|f\|_{C(D)}=\|a\|_{\ell_2(\R^n)}$, 
 the map $a\rightarrow f_a(x)$ establishes a linear isometry between $V$ with the supremum norm and 
 $\R^n$ with the Euclidean norm. 
 Let   $M$ be the measurement map given by  the linear functionals corresponding to point evaluation at any set  $\{P_j\}_{j=1}^m$  of $m$ points from $D$. Then $M$ maps $C(D)$ into 
$\ell_\infty(\R^m)$ and $\|M_V\|=1$.
It follows from \eref{tu}  that   $\mu(\cN,V)\approx \|M_V\|\cdot \|M_V^{-1}\|$. 
This means that  
$$
\mu(\cN,V)\leq C d(\ell_2(\R^n), M(V)), \quad M(V)\subset \ell_\infty(\R^m),
$$
where $d(\ell_2(\R^n), M(V)):=\inf \{\|T\|\|T^{-1}\|, \, T:\ell_2(\R^n)\rightarrow M(V), \, T\,\, \mbox{isomorphism}\}$ 
is the  Banach-Mazur distance between  the  $n$ dimensional Euclidean space $\R^n$ and  the subspace 
$M(V)\subset \ell_\infty(\R^m)$. 
 It  is a well known, but nontrivial  fact  in the local theory of Banach spaces (see \cite[Example 3.1]{FLM} 
 or \cite[Section 5.7]{MSch}) that  to keep 
  $d(\ell_2(\R^n), M(V))\leq C$, one needs   $\ln m\geq c n$.

The scenario of the last paragraph is  the worst that can happen.   To see why, let us 
recall the following notion:  a set  $A$ is a $\delta$-net for a set  $S$ 
($A\subset S\subset \cX$ and $\delta>0$)   if for every $x\in S$ there exists $a\in A$, such that $\|x-a\|\leq \delta$.
For a given $n$-dimensional subspace $V\subset \cX$ and $\delta>0$, let us fix a 
$\delta$-net  $\{v_j\}_{j=1}^N$  for  $ \{v\in V\,:\,\|v\|=1\}$ with $N\leq (1+2/\delta)^n$. It is well known that 
such a net exists (see  \cite[Lemma 2.4]{FLM} or \cite[Lemma 2.6]{MSch}).  Let $l_j\in \cX^*$   
be  norm one functionals,  such that $1=l_j(v_j)$, $j=1,2,\dots,N$. 
We define our measurement $M$ as $M=(l_1,\ldots,l_N)$, so $\cN=\bigcap_{j=1}^N \ker l_j$. 
When $x\in \cN$, $v\in V$  with $\|v\|=1$,  and $v_j$ is such that $\|v-v_j\|\leq \delta$, we have
$$
\|x-v\|\geq \|x-v_j\|-\delta\geq |l_j(x-v_j)|-\delta =1-\delta,
$$
and so  for this choice of $M$, we have
$$
 \mu(\cN,V)\leq 2\mu(V,\cN)=2 \sup_{x\in \cN,\,  v\in V,\,\|v\|=1\,}\frac{1}{\|x-v\|}\leq \frac{2}{1-\delta}.
$$
\begin{remark}
\label{rem:betterm}For specific Banach spaces $\cX$ and subspaces $V\subset \cX$, the situation is much better.   
We have already discussed 
such example in the case of the space of trigonometric
polynomials and $\cX$ the space of periodic functions in $C([-\pi,\pi])$.
\end{remark}


\begin{remark}  
Let us discuss
briefly the situation when  again $\cX=C([-\pi,\pi])$,  but now the measurements are 
given as lacunary Fourier coefficients,  i.e.,  
 $M(f)=(\hat f(2^0), \ldots,\hat f(2^m))$.
 From \eref{Quotdual2}, we infer that
$\|\alpha\|_M^*=\frac1{2\pi}
\int_{-\pi}^{\pi}|\sum_{j=0}^m \alpha_j e^{i2^j
t}|\, dt$. Using a well known analog of Kchintchine's
inequality valid for lacunary trigonometric polynomials (see,  e.g. {\rm\cite[ch 5, Th.8.20]{Z})},  we
derive
$$
\sum_{j=0}^m|w_j|^2\leq \|w\|_M^2\leq C\sum_{j=0}^m|w_j|^2,
$$
for some constant $C$. If $\Delta:\R^{m+1}\rightarrow C([-\pi,\pi])$ is any linear lifting, using 
{\rm\cite[III.B.5 and III.B.16]{Wbook}}, we obtain
$$
\|\Delta\|=\|\Delta\|\cdot\|M\|\geq \frac{\sqrt\pi}{2C}\sqrt{m+1}.
$$
On the other hand,
there exists a constructive, nonlinear lifting
$\Delta_F:\R^{m+1}\rightarrow C([-\pi,\pi])$
with $\|\Delta_F\|\leq \sqrt{e}$,
see {\rm \cite{F}}.
\end{remark} 
 


\section{Appendix: Improved estimates for $\rad(\cK_w)$}

The purpose of this appendix is to show that the estimates derived in Theorem \ref{theorem:main} can be 
improved if we assume more structure for the
Banach space $\cX$.  Let us  begin by recalling a lemma from  \cite[Lemma 4]{Figiel} (see also \cite[Cor. 2.3.11]{AOS}).

\begin{lemma}
\label{FL4} 
If $x,y\in \cX$ and $\|x\|,\|y\|\leq \e$,  then  $\|x+y\|\leq 2\e [1-\delta_\cX(\|x-y\|/\e)]$ whenever $\e>0$.
\end{lemma}

\noindent
As noted in the introduction, when  $\cX$ is uniformly convex, then  $\delta_\cX$ is strictly increasing.  Therefore,  it has a compositional inverse $\delta_\cX^{-1}$ which is also strictly increasing and satisfies $\delta^{-1}_\cX(t)\le 2$, $0\le t\le 1$  and $\delta^{-1}_\cX(0)=0$.  Hence, the   following result improves on the upper estimate (ii) of Theorem \ref{theorem:main}.

\begin{prop} 
\label{USprop2}
Let $\cX$   be a uniformly convex Banach space with modulus of convexity $\delta_\cX$, defined in \eref{uc}.
If 
$$
\min_{x\in \cK_w}\dist(x,V)=\gamma\e, \,\, \mbox{  for  } \gamma\leq1,
$$
 then, 
$$
 \diam(\cK_w)\leq  \e \mu(\cN,V) \delta_{\cX}^{-1}(1-\gamma).
$$

\end{prop}
{\bf Proof:} Let us fix any $u_1,u_2\in \cK_w$ and take $u_0=\frac{u_1+u_2}{2}$.
For any $u\in \cX$,  we denote by $P_V(u)$ the best approximation to $u$ from $V$,  which is unique.  
Then, we have
\begin{eqnarray}
\label{firstestimate}
\dist(u_0,V)&=&\|u_0-P_{V}(u_0)\|\leq \|\frac{u_1+u_2}{2}-\tfrac12 (P_{V}(u_1)+P_{V}(u_2)) \|\nonumber\\
&=&\left\|\tfrac12  ([u_1-P_{V}(u_1)]+[u_2-P_{V}(u_2)]) \right\|:=\tfrac12\|\eta_1 +\eta_2\|.
\end{eqnarray}
  Let $\alpha :=\|\eta_1-\eta_2\|$.
Since   $\|\eta_1\|,\|\eta_2\|\leq \e$, if we start with \eref{firstestimate} and then using Lemma \ref{FL4},   we find that 
\begin{eqnarray*}
\gamma\e \leq \dist(u_0,V) \leq \tfrac12\|\eta_1+\eta_2\| \leq  \e [1-\delta_\cX(\alpha/\e)].
\end{eqnarray*}
This gives  $\gamma\leq 1-\delta_\cX(\alpha/\e)$, so  we get 
$$
\alpha\leq \e\delta_\cX^{-1}(1- \gamma).
$$
Clearly $u_1-u_2\in\cN$,  so we have
$[P_V(u_1)-P_V(u_2)]+[\eta_1-\eta_2]=u_1-u_2\in \cN$.
From the definition of $\mu(\cN,V)$, see \eref{defmu}, we derive that
$$
\mu(\cN,V)\geq  \frac{\|u_1-u_2\|}{\|u_1-u_2-P_V(u_1)+P_V(u_2)\|}=  \frac{\|u_1-u_2\|}{\|\eta_1-\eta_2\|},
$$
which gives
$$
\|u_1-u_2\|\leq \alpha \mu(\cN,V)\leq \e \mu(\cN,V) \delta_{\cX}^{-1}(1-\gamma) .
$$
Since $u_1,u_2$ were arbitrary we have proven the claim.
 \hfill $\blacksquare$

We next give an estimate of $\diam(\cK_w)$ from below by using the concept of the   modulus of  smoothness.   Recall that 
the modulus of smoothness of   $\cX$ (see   e.g. \cite{Figiel,LT,AOS})  
is defined by 
\begin{equation}\label{USdef}
\rho_\cX(\tau)=\sup\left\{\frac{\|x+y\|+\|x-y\|}{2}-1\ :\ \|x\|=1,\, \|y\|\leq \tau\right\}.
\end{equation}
The space $\cX$ is said to be {\it uniformly smooth} if $\displaystyle{\lim_{\tau\to 0}\rho_\cX(\tau)/\tau =0}$. 
Clearly,   $\rho_\cX(\tau)$ is a pointwise supremum of a family of convex functions, so it is a convex and  strictly  increasing function of $\tau$ (see \cite[Prop. 2.7.2]{AOS}).
Let us   consider the quotient space $\cX/V$. Each element of this space is a coset $[x+V]$ with the norm
$\|[x+V]\|_{\cX/V}:=\dist(x,V)$. It is known  (for example,  it easily follows from  \cite[Prop.1.e.2-3]{LT}) 
that $\rho_{\cX/V}(\tau)\leq \rho_\cX(\tau)$, 
so any quotient space of a uniformly smooth space is also uniformly smooth.

\begin{lemma}\label{US2}
Suppose  $\cX$ is  a  Banach space.  If $u_0,u_1\in \cX$ with  $\|u_0\|=\|u_1\|=\e$ and
 $$
 \inf_{\lambda\in(0,1)}\|\lambda u_0+(1-\lambda)u_1\|=\gamma \e,\quad 0<\gamma\leq 1,
 $$
  then
\begin{equation}
\label{US1} 
\|u_0-u_1\|\geq 2 \gamma \e \rho_\cX^{-1}(\tfrac{1-\gamma}{2\gamma}).
\end{equation}
\end{lemma}

\noindent
{\bf Proof:} Let  $x:=\lambda_0 u_0+(1-\lambda_0) u_1$  be such that $\|x\|=\gamma\e$. We 
 denote by $\bar x:=x/\gamma\e$,  $y:=\frac{u_0-u_1}{2 \gamma\e}$ and  let $\tau=\|y\|$. 
It follows from   (\ref{USdef}) that
\begin{equation}\label{US3a}
\rho_\cX(\tau)\geq \left(\frac{\|\bar x+y\|+\|\bar x-y\|}{2}-1\right)=\frac{1}{2\gamma\e}(\|x+\tfrac12(u_0-u_1)\|+\|x-\tfrac12(u_0-u_1)\|)-1. 
\end{equation}
Consider the  function $\phi(t)=\|x+t(u_0-u_1)\|:=\|x_t\|$ defined for $t\in \R$. It is a convex function which attains its minimum 
value  $\gamma\e$ when  $t=0$. For $t=-\lambda_0$ we get $x_t=u_1$ and for $t=1-\lambda_0$ we get $ x_t=u_0$. This implies that for $t\in[-\lambda_0,1-\lambda_0]$ we have $\gamma\e\leq \phi(t)\leq \e$ and for $t\notin[-\lambda_0,1-\lambda_0]$ we have $\phi(t)\geq \e$. So from (\ref{US3a}),  we get
$$
\rho_\cX(\tau)\geq\frac{1}{2\gamma\e}  (\gamma\e+\e)-1=\frac{1-\gamma}{2\gamma}, 
$$
 which yields (\ref{US1}). \hfill $\blacksquare$

With the above lemma in hand, we can give the following lower estimate for the diameter of $\cK_w$.
\begin{prop}\label{USprop1}
Let $\cX$ be a Banach space  with modulus of smothness $\rho_\cX$.
If
$$
\min_{x\in \cK_w}\dist(x,V)=\gamma\e,\,\, \mbox{  for  }\quad 0< \gamma\leq1,
$$
 then
 $$
\diam (\cK_w)\geq 2 \e \mu(\cN,V) \gamma\rho_\cX^{-1}(\tfrac{1-\gamma}{2\gamma}).
$$
 \end{prop}
 
 \noindent
{\bf Proof:}
Let  $x_0 \in \cK_w$ satisfy $\dist(x_0,V)=\gamma\e$.  Given any  $\eta\in \cN$ with $\|\eta\|=1$,
we fix $\alpha,\beta>0$,  such that $\dist(x_0+\alpha \eta,V)=\e=\dist(x_0-\beta \eta,V)$.
Note that $x_0+\alpha \eta$ and $x_0-\beta \eta$ belong to $\cK_w$, and therefore,
\be
\label{USa3}
 \diam (\cK_w) \ge \alpha+\beta.
\ee

\noindent
We now  apply Lemma \ref{US2} for  the quotient space $\cX/V$ with $u_0=[x_0+\alpha \eta+V]$ 
and $u_1=[x_0-\beta \eta+V]$.   It follows from \eref{US1}  that 
\be\label{USa1}
\dist((\alpha+\beta)\eta,V)=\|u_0-u_1\|_{\cX/V}\geq 2\e \gamma \rho_{\cX/V}^{-1}(\tfrac {1-\gamma}{2\gamma})
\geq 2\e \gamma \rho_{\cX}^{-1}(\tfrac {1-\gamma}{2\gamma}).
\ee
Finally, observe that, in view of \eref{USa3}, we have 
\begin{eqnarray*}
\inf_{\eta \in \cN,\, \|\eta\|=1}(\alpha+\beta)\dist (\eta,V)\leq  \diam(\cK_w)  \inf_{\eta\in \cN,\, \|\eta\|=1}\dist(\eta,V)
= \diam(\cK_w)\Theta(\cN,V) = \frac{\diam(\cK_w)}{\mu(\cN,V)}. 
\end{eqnarray*}
Therefore,
 using \eref{USa1}, we arrive at 
$$
\diam (\cK_w)\geq  2\e \mu(\cN,V) \gamma\rho_\cX^{-1}(\tfrac{1-\gamma}{2\gamma}).
$$
 \hfill $\blacksquare$

\begin{remark}  For a general Banach space  $\cX$,    we have that $\delta_\cX(\tau)\geq 0 
$ and 
$\rho_\cX(\tau)\leq \tau$, for $\tau>0$,
and in general those are the best estimates. 
  So for every Banach space $\cX$ we obtain from Proposition {\rm \ref{USprop1}} that 
$\diam(\cK_w)\geq \e(1-\gamma)\mu(\cN,V)$.

Moduli of convexity and smoothness are computed (or well estimated)  for various classical spaces.
In particular, their exact values  for the $L_p$ spaces, $1<p<\infty$, have been  computed in {\rm\cite{Hanner}}.     
We will just state the asymptotic results (see e.g. {\rm\cite{LT}}) 
\begin{eqnarray*}
\delta_{L_p}(\e)&=& \begin{cases}(p-1)\e^2/8 +o(\e^2), &\mbox{ for  } 1<p<2,\\
\e^p/p2^p +o(\e^p),& \mbox{ for  } 2\leq p<\infty,
\end{cases}\\
\rho_{L_p}(\tau)&=& \begin{cases}\tau^p/p+o(\tau^p), &\mbox{  for  } 1<p\leq 2,\\
(p-1)\tau^2/2+o(\tau^2), & \mbox{  for  } 2\leq p <\infty.
\end{cases}
\end{eqnarray*}
From the parallelogram identity, we have,
$$
\delta_{L_2}(\e)=1-\sqrt{1-\e^2/4},\quad \rho_{L_2}(\tau)=\sqrt{1+\tau^2}-1.
$$  
 It follows from Propositions  {\rm \ref{USprop2}} and {\rm\ref{USprop1}} that
 $$
\e \mu(\cN,V)\sqrt{1-\gamma^2}\leq \diam(\cK_w)_{L_2}\leq2\e \mu(\cN,V)\sqrt{1-\gamma^2}.
 $$
By isomorphism of Hilbert spaces this last result holds for any Hilbert space, and (up to a constant) we retrieve the results of 
{\rm\cite{BCDDPW1}}.
\end{remark}

 \noindent
Ronald DeVore\\
Department of Mathematics, Texas A\&M University,
College Station, TX 77840, USA\\
  rdevore@math.tamu.edu

 \vskip .1in
\noindent
Guergana Petrova\\
Department of Mathematics, Texas A\&M University,
College Station, TX 77840, USA\\
gpetrova @math.tamu.edu.
 \vskip .1in
\noindent
Przemyslaw Wojtaszczyk\\
Interdisciplinary Center for Mathematical and Computational Modelling, \\
University of Warsaw, 00-838 Warsaw, ul. Prosta 69, Poland\\
 wojtaszczyk@icm.edu.pl
\\

  \end{document}